\documentclass{article}

\usepackage{arxiv}

\usepackage[utf8]{inputenc}
\usepackage[T1]{fontenc}
\usepackage[english]{babel}

\usepackage{lineno}
\usepackage{standalone}
\usepackage{caption}
\usepackage{subcaption}
\usepackage{stmaryrd}
\usepackage{amssymb,amsmath,amstext,empheq}
\usepackage{upgreek}
\usepackage{graphicx}
\usepackage{booktabs}

\usepackage{tikz}
\usetikzlibrary{decorations.pathmorphing,decorations.pathreplacing,decorations.markings,patterns,calc}
\usetikzlibrary{shapes.geometric, shapes.arrows, arrows, angles, quotes}
\usetikzlibrary{spy,backgrounds}
\usetikzlibrary{external}
\usetikzlibrary{intersections}
\usetikzlibrary{tikzmark}

\tikzstyle{startstop} = [rectangle, rounded corners, minimum width=3cm, minimum height=1cm,text centered, draw=black, fill=red!30]
\tikzstyle{io} = [trapezium, trapezium left angle=70, trapezium right angle=110, minimum width=3cm, minimum height=1cm, text centered, draw=black, fill=blue!30]
\tikzstyle{process} = [rectangle, minimum width=3cm, minimum height=1cm, text centered, draw=black, fill=orange!30]
\tikzstyle{decision} = [diamond, minimum width=3cm, minimum height=1cm, text centered, draw=black, fill=green!30]
\tikzstyle{arrow} = [->,>=stealth]

\usepackage{import}
\usepackage{pdfpages}
\usepackage{transparent}
\usepackage{xcolor}

\usepackage[normalem]{ulem}
\usepackage{textcomp}
\usepackage{siunitx}
\usepackage[most]{tcolorbox}

\usepackage{bm}
\usepackage[]{datatool}
\usepackage{calc}
\usepackage{textgreek}
\usepackage{bigints}
\usepackage{physics}
\usepackage{interval}
\usepackage{accents}
\usepackage{mathtools,eqparbox}
\usepackage{indentfirst}
\usepackage[algo2e,ruled,linesnumbered,vlined,commentsnumbered,longend]{algorithm2e}
\SetKw{Continue}{continue}
\SetKw{Break}{break}

\usepackage{enumerate}

\usepackage{relsize}
\usepackage[column=O]{cellspace}
\usepackage{pgfplots}
\pgfplotsset{compat=newest}
\usepgfplotslibrary{groupplots,dateplot}
\usepgfplotslibrary{colormaps}
\usepgfplotslibrary{colorbrewer}
\usepgfplotslibrary{polar}
\usepgfplotslibrary{fillbetween}

\usepackage[absolute,overlay]{textpos}

\usepackage{tikz-dimline}

\usepackage[autostyle]{csquotes}

\usepackage[colorinlistoftodos]{todonotes}

\usepackage{hyperref}
\usepackage{hypcap}
\usepackage{bookmark}

\DeclareMathOperator*{\argmin}{arg\,min}

\DeclareMathOperator*{\argmax}{arg\,max}

\DeclareMathOperator*{\Div}{div}

\DeclareMathOperator{\sgn}{sgn}

\DeclareMathOperator{\Grad}{grad}

\newcommand{\bu}{\boldsymbol{u}} 
\newcommand{\bn}{\boldsymbol{n}} 
\newcommand{\bx}{\boldsymbol{x}} 

\newcommand{\btau}{\boldsymbol{\tau}} 

\newcommand{\ba}{\boldsymbol{a}}

\newcommand{\bsigma}{\boldsymbol{\sigma}} 
\newcommand{\bepsilon}{\boldsymbol{\varepsilon}} 

\newcommand{\bb}{\boldsymbol{b}}

\newcommand{\jump}[1]{\left\llbracket #1 \right\rrbracket} 
\newcommand{\bphi}{\boldsymbol{\varphi}}

\newcommand{\rom}[1]{\uppercase\expandafter{\romannumeral #1\relax}}

\intervalconfig {
	soft open fences ,
}

\usepackage{amsthm}

\theoremstyle{plain}

\theoremstyle{definition}
\newtheorem{definition}{Definition}[section]

\theoremstyle{remark}

\pgfplotsset{
	legend image with text/.style={
			legend image code/.code={%
					\node[anchor=center] at (0.3cm,0cm) {#1};
				}
		},
}

\usepackage{doi}
\usepackage[
	natbib=true,
	backend=biber,
	giveninits=true,
	date=year,
	doi=false,
	url=false,
	isbn=false,
	eprint=false,
	bibstyle=ieee,
	citestyle=numeric-comp,
]{biblatex}
\title{
	On optimization of heterogeneous materials for enhanced resistance to bulk fracture
}
\author{
  Sukhminder Singh$^{1, 2}$\thanks{Corresponding author, email address: \texttt{sukhminder.singh@fau.de}}\ , Lukas Pflug$^1$, Julia Mergheim$^3$, and Michael Stingl$^2$ \\ \\
$^1$Competence Unit for Scientific Computing (CSC), \\
Friedrich-Alexander-Universit\"at~Erlangen-N\"urnberg~(FAU) \\ \\
$^2$Chair of Applied Mathematics (Continuous Optimization), \\
Friedrich-Alexander-Universit\"at~Erlangen-N\"urnberg~(FAU) \\ \\
$^3$Institute of Applied Mechanics, \\
Friedrich-Alexander-Universit\"at~Erlangen-N\"urnberg~(FAU)
}

\begin{document}









\maketitle
\begin{abstract}
	We propose a novel approach to optimize the design of heterogeneous materials, with the goal of enhancing their effective fracture toughness under mode-\rom{1} loading.
	The method employs a Gaussian processes-based Bayesian optimization framework to determine the optimal shapes and locations of stiff elliptical inclusions within a periodic microstructure in two dimensions.
	To model crack propagation, the phase-field fracture method with an efficient interior-point monolithic solver and adaptive mesh refinement, is used.
	To account for the high sensitivity of fracture properties to initial crack location with respect to heterogeneities, we consider multiple cases of initial crack and optimize the material for the worst-case scenario.
	We also impose a minimum clearance constraint between the inclusions to ensure design feasibility.
  Numerical experiments demonstrate that the method significantly improves the fracture toughness of the material compared to the homogeneous case.
\end{abstract}


\section{Introduction}\label{sec:introduction}

The field of mathematical and computational fracture mechanics has seen significant advancements in the past two decades, leading to a comprehensive understanding of the fracture phenomenon in engineering structures at various scales. These advancements have facilitated the exploration of several mechanisms that govern crack nucleation, growth, and ultimate structural failure.
Sharp cracks can form during the manufacturing process or arise during service at stress concentration regions, owing to several factors, such as structural overloads, fatigue, and corrosion. When the stresses at the crack tip exceed the critical stress defined by the material's fracture toughness, the crack can propagate through the structure, potentially leading to unforeseen structural collapse.

Introducing architected heterogeneities or voids in materials at the microscale is a promising approach to enhance the fracture resistance and load-carrying capacity of structures.
Such heterogeneities can effectively retard or even arrest crack propagation~\cite{TowardsBrittleLebiha2021,EffectiveToughHossai2014}.
The size, shape, and material properties of the heterogeneities can be chosen based on the solution of a numerical optimization~\cite{NumericalOptim2006} problem, in which the objective is to maximize or minimize a quantity of interest that characterizes the fracture toughness of the heterogeneous material.

However, the quasi-static analysis of crack propagation exhibits two significant characteristics, namely snap-back/snap-through instabilities and bifurcations~\cite{wriggers2008nonlinear,NonlinearDynamStroga2018}.
These behaviors are also observed in other structural problems such as nonlinear buckling~\cite{StabilityBifuPignat1998}, crash~\cite{SurrogateModelBoursi2018}, and damage~\cite{StabilityOfStBazant2010}.
In such cases, the required derivatives are not readily available, making the use of standard gradient-based optimization methods infeasible.
The presence of numerical noise in fracture simulations further poses a challenge for the optimization algorithms.
Therefore, developing derivative-free and robust optimization schemes to handle such complex and nonlinear problems is crucial for designing materials with architected heterogeneities.

Efforts have been made to simplify the structural analysis model, and use density- (e.g., SIMP~\cite{GeneratingOptiBendso1988}) and level-set-based~\cite{ALevelSetMetAllair2002} topology optimization techniques to obtain stronger structures that are more resistant to fracture.
Stress constraints have been incorporated historically to achieve success in improving resistance to crack nucleation in areas of high stress concentration~\cite{StudyOnTopoloCheng1992,TopOptStressConstraintsBendsoe1998,StressBasedToLeCh2010,AnEnhancedAggLuoY2013,StressConstraiAmir2017,StressBasedShPicell2018}.
Similar techniques have been developed for stationary cracks to improve their resistance to grow further~\cite{OptimizationOfGuGr2016,TopologyOptimiKang2017,MinimizingCracKlarbr2018,FractureStrengHuJi2020,OnTailoringFrZhang2022}.
For such problems, efficient gradient-based optimization methods can reliably exploit the objective's analytical sensitivities concerning the design parameters.


In recent years, there have been numerous research efforts aimed at integrating propagating cracks into optimization frameworks.
The primary approach involves regularizing the description of discrete cracks through the use of a scalar phase field
~\cite{RevisitingBritFrancf1998,AReviewOnPhaAmbati2015}.
Typically, the objective function consists of a structural response function integrated over a finite number of fixed load steps, such as
integrated mechanical work~\cite{TopologyOptimiDaDa2018,TopologyOptimiDaDa2020,LevelSetTopolWuCh2020,APathDependenWuCh2021},
integrated fracture surface energy~\cite{ANovelTopologRuss2020}, or
integrated elastic energy~\cite{TopologyOptimiDesai2022}.
While these methods have demonstrated the ability to enhance the fracture resistance of structures, they necessitate some regularization of the fracture problem and do not account the handling of unstable crack propagation.

On the other hand, the fracture toughness criteria for heterogeneous materials, as proposed by Hossain~et~al.~\cite{EffectiveToughHossai2014}, involves calculating the J-integral,
whose optimization is still an unexplored area.
The method involves recording the J-integral curve with respect to mode-\rom{1} surfing boundary conditions and defining the effective fracture toughness as the peak of the J-integral curve.
However, optimizing the peak J-integral poses several challenges.
First, the objective function is not differentiable with respect to design parameters due to the presence of unstable crack propagation.
Second, the optimization algorithm should consider multiple crack patterns dependent on the choice of the location of the initial crack.
Third, the fracture problem itself should be computationally tractable given the high computational costs of solving the phase-field fracture problem.

In this work, we propose a Gaussian processes-based Bayesian optimization strategy to optimize the shapes and placement of stiff elliptical inclusions within a periodic microstructure in 2D, aiming to maximize its effective fracture toughness~\cite{EffectiveToughHossai2014} with respect to matrix cracking under mode-I loading.
The phase-field fracture method is utilized, while the interfaces are assumed to be linear elastic.
The extended finite element method is used for discretization of the evolving interfaces, which makes it suitable for use in a parametric framework without the need for finite element remeshing.
To reduce the computational burden of the fracture simulation, an interior-point monolithic solver is used instead of the standard alternating minimization algorithm~\cite{NumericalExperBourdi2000,TheVariationalBourdi2008,NumericalImpleBourdi2007,ValidationSimuMesgar2015,AnAdaptiveFinBurke2010} for phase-field fracture, along with adaptive mesh refinement~\cite{APrimalDualAHeiste2015,AnAdaptiveMesGupta2022} near the crack tip and subsequent coarsening along the tail of the crack.
To ensure design feasibility, an inequality constraint is introduced, which requires a minimum clearance between the inclusions.
In addition, the optimization algorithm considers multiple initial crack locations to account for the possibility of different crack patterns~\cite{StochasticPhasGerasi2020,DeterministicANagara2023} realized for a given design, and optimizes for the worst-case of the corresponding objective values.

The remainder of the paper is structured as follows.
Section~\ref{sec:fracture_toughness} gives a concise definition of the effective fracture toughness of heterogeneous materials, as proposed by Hossain~et~al.~\cite{EffectiveToughHossai2014}.
In Section~\ref{sec:fracture_problem_formulation}, we formulate a phase-field fracture model based on energy minimization principles.
Section~\ref{sec:spatial_discretization} offers a brief overview of the spatial discretization scheme, including the extended finite element method used to model sharp interfaces.
Section~\ref{sec:simulation_algorithm} presents the solution algorithm for the phase-field fracture problem.
In Section~\ref{sec:optimization}, the design optimization problem formulation is presented, including design parametric and optimization studies.
Finally, Section~\ref{sec:conclusion} concludes the paper.


\section{Mode-\rom{1} fracture toughness of heterogeneous materials}\label{sec:fracture_toughness}


\begin{figure}[ht]
	\centering
	\def\svgwidth{1\columnwidth}
	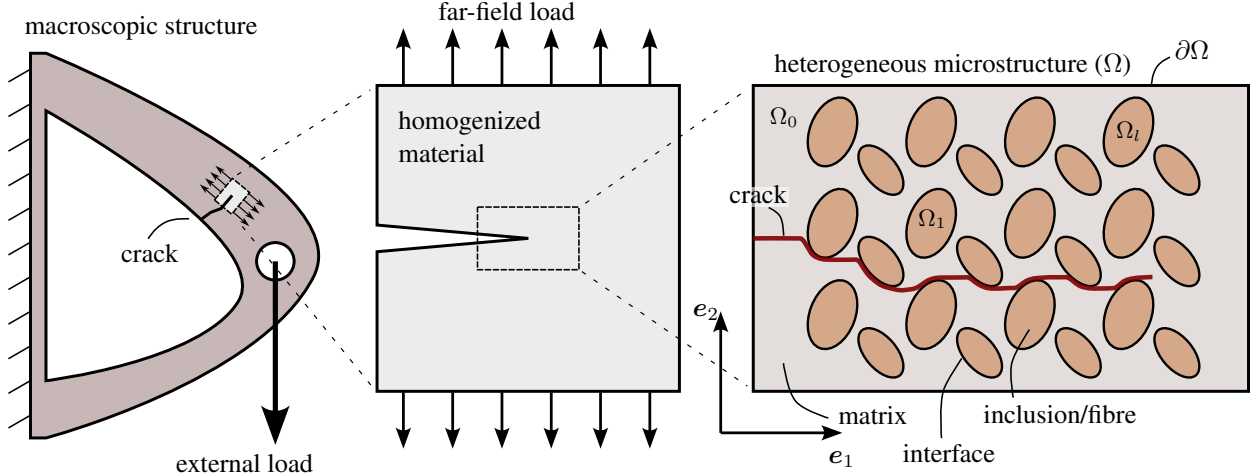

	\caption{
		Illustration of a crack under mode-\rom{1} loading condition at the mesoscale (middle). At the macroscale, the crack grows through crack propagation mechanisms (fatigue, corrosion, etc.), while at the microscale, the crack interacts with the heterogeneities that provide resistance to its propagation.
		(Figure adapted from Hansen-D\"orr~\cite{doerr2022phase}).
	}
	\label{fig:multiscale_schematic}
\end{figure}

In order to calculate the mode-\rom{1} fracture toughness of heterogeneous materials, Hossain~et~al.~\cite{EffectiveToughHossai2014} proposed a methodology that involves subjecting the microstructure of the material, represented by the domain $\Omega$ shown in Figure~\ref{fig:multiscale_schematic}, to surfing boundary conditions applied to the domain boundary $\partial\Omega$.
For crack propagation at the macroscale, these boundary conditions are defined by the displacement field:
\begin{subequations}
	\begin{align}
		\hat{u}_1(\bx, t^{(k)}) & \coloneqq \frac{K_\mathrm{\rom{1}}}{2\mu} \sqrt{\frac{r(\bx, t^{(k)})}{2\pi}} \left[\kappa - \cos(\theta(\bx, t^{(k)})) \right] \cos(\frac{\theta(\bx, t^{(k)})}{2}), \\
		\hat{u}_2(\bx, t^{(k)}) & \coloneqq \frac{K_\mathrm{\rom{1}}}{2\mu} \sqrt{\frac{r(\bx, t^{(k)})}{2\pi}} \left[\kappa - \cos(\theta(\bx, t^{(k)})) \right] \sin(\frac{\theta(\bx, t^{(k)})}{2}),
	\end{align}
	\label{eq:surfing_boundary_conditions}
\end{subequations}
where $\bx$ and $t^{(k)}$ are the spatial and time variables, respectively, with $k$ representing the load or pseudo-time step.
The material parameter $K_\mathrm{\rom{1}}$ is the mode-\rom{1} stress-intensity factor,
$\kappa=3-4\nu$ is an elastic constant for plane-strain case~\cite{FiniteElementsKuna2013}, $\mu=E/2(1+\nu)$ is the shear modulus, $E$ is the elastic modulus and $\nu$ is the Poisson's ratio.
These material properties are chosen to be representative of the matrix.

The polar coordinates $(r, \theta)$, which emerge from the crack tip, are expressed as
\begin{subequations}
	\begin{align}
		r(\bx, t^{(k)})      & \coloneqq \sqrt{{\left(x_1 - vt^{(k)}\right)}^2 + \left(x_2 - \bar{x}_2^{(k-1)}\right)^2},          \\
		\theta(\bx, t^{(k)}) & \coloneqq \tan^{-1}\left(\frac{x_2 - \bar{x}_2^{(k-1)} }{x_1 - vt^{(k)}}\right),
	\end{align}
	\label{eq:polar_coordinates_surfing}
\end{subequations}
where $v$ denotes the macroscopic crack velocity,
and $\bar{\bx}^{(k)}$ represents the position of the crack tip, with $\bar{\bx}^{(0)} = \boldsymbol{0}$.

Given the mode-\rom{1} crack loading conditions, the J-integral~\cite{APathIndependRice1968,CrackPropagatiCherep1967} can be utilized to determine the crack driving force in heterogeneous media.
Specifically, the J-integral is calculated along the boundary of the domain using the following expression:
\begin{align}
	\mathcal{J}(\bu) \coloneqq \int_{\partial \Omega} \boldsymbol{e}_1 \cdot \boldsymbol{C}(\bu(\bx)) \cdot \boldsymbol{n}(\bx) \, \dd{s},
\end{align}
where
\begin{align}
	\boldsymbol{C}(\bu) \coloneqq \Psi(\bu) \boldsymbol{I} - \nabla^t \boldsymbol{u} \cdot \bsigma(\bu)
\end{align}
represents Eshelby's energy-momentum tensor~\cite{TheElasticEneEshelb1975} or the configurational stress tensor, $\Psi(\bu) = \frac{1}{2} \bsigma(\bepsilon(\bu)) :\boldsymbol{\varepsilon}(\bu)$ denotes the elastic energy density, $\bsigma$ is the second-order stress tensor and $\boldsymbol{\varepsilon}$ correspond to the second-order strain tensor.
The effective fracture toughness $G_\mathrm{c}^\text{eff}$ of the heterogeneous material is then defined as the maximum value of the J-integral observed over time under the surfing boundary conditions, i.e.,
\begin{align}
  G_\mathrm{c}^\text{eff} \coloneqq \max_{k} \mathcal{J}(\bu^{(k)}).
  \label{eq:effective_fracture_toughness}
\end{align}
%

\section{Phase-field fracture problem formulation}\label{sec:fracture_problem_formulation}

\begin{figure}[ht]
	\centering
	\def\svgwidth{0.8\columnwidth}
	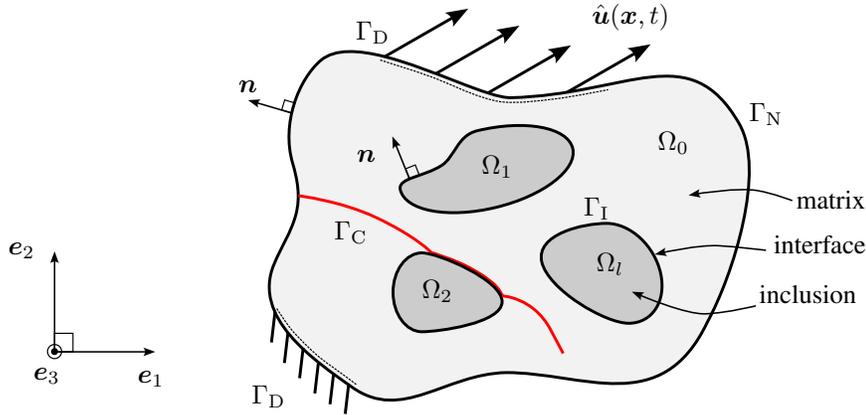

	\caption{Schematic of heterogeneous continuum domain with crack $\Gamma_\mathrm{C}$.}
	\label{fig:continuum_domain}
\end{figure}

Let $\Omega \subset \mathbb{R}^d, d=\{2, 3\}$ be a smooth, open, connected and bounded domain representing the configuration of a composite body, as shown in Figure~\ref{fig:continuum_domain}.
The body consists of a matrix sub-domain $\Omega_0$ and $L \in \mathbb{N}$ disjoint subsets $\Omega_l, l = 1,\dots,L$ representing inclusions embedded in the matrix.
The inclusions are bonded to the matrix material by elastic interfaces denoted by the set $\Gamma_\mathrm{I} \coloneqq \cup_{l=1}^L \partial \Omega_l$.
The domain boundary $\partial \Omega$ is divided into two disjoint sets, $\Gamma_\mathrm{D}$ and $\Gamma_\mathrm{N}$, representing Dirichlet and Neumann boundaries, respectively.
Let us assume that $\Omega$ contains the crack $\Gamma_\mathrm{C}$, which is a $d-1$ dimensional, possibly disconnected, set.

The state of this mechanical system for a given set of boundary conditions is given by the displacement field $\bu : \overline{\Omega} \rightarrow \mathbb{R}^d$, which admits discontinuities on $\Gamma_\mathrm{C}$ and $\Gamma_\mathrm{I}$.
The interfacial opening at a point $\bx \in \Gamma_\mathrm{I}$ with unit vector $\bn(\bx)$ pointing in the direction perpendicular to the interface is given by the displacement jump $\jump{\bu(\bx)} \coloneqq \bu(\bx^+) - \bu(\bx^-)$, where
$\bx^\pm = \lim_{\epsilon \rightarrow 0^+} \bx \pm \epsilon \bn(\bx)$.
Furthermore, we assume that $\Gamma_\mathrm{I}^+$ denotes the matrix side of the interface and $\Gamma_\mathrm{I}^-$ denotes the inclusion side of the interface.
Also, we assume that the normal vector $\bn$ always points in the direction away from the inclusion.

On the Dirichlet boundary, a time-dependent displacement field $\hat{\bu} : \Gamma_\mathrm{D} \times [0, T] \rightarrow \mathbb{R}^d, \ T \in \mathbb{R}_{> 0}$ is prescribed, whereas for the sake of simplicity of the model, we assume that no external body and surface forces act on the system.

The phase-field fracture model is founded on the idea that the fracture process can be described by a continuous field defined over the full material domain, known as the phase field or \emph{damage} field, denoted by $\alpha : \Omega \rightarrow \mathbb{R}$.
In the uncracked region, the phase field is \num{0} and in the fully cracked zone, it is \num{1}.
The damage is characterized by deterioration of the bulk stiffness which approaches \num{0} at the points where $\alpha$ approaches \num{1}.
We define function spaces for the displacement and damage fields as
\begin{align}
	\mathcal{S} \coloneqq \left\{ \bu \in L^2(\Omega; \mathbb{R}^d) : \bu \vert_{\Omega_l} \in H^1(\Omega_l; \mathbb{R}^d),\ l=1,\dots,L \right\},
	\label{eq:state_solution_space}
\end{align}
and
\begin{align}
	\mathcal{A} \coloneqq \left\{ \alpha \in H^1(\Omega; \mathbb{R}) : \alpha \in [0, 1], \alpha = 0 \text{ on } \partial \Omega \right\},
\end{align}
respectively.
By restricting the value of $\alpha$ to zero at the domain boundary, the formation of cracks at the boundary is prevented as it requires only half of the fracture energy for a crack to propagate along the boundary.
To simplify the model, it is assumed that the stiffness within the inclusion domain $\Omega_l, l=1,\dots,L$ does not degrade, and hence, cracks do not appear within the inclusions. Furthermore, it is assumed that the fracture toughness of both the interfaces and the bulk material are equivalent. Consequently, the interfacial cracks can be represented using the same damage field used for the bulk cracks, albeit in a diffused manner.

Following Bourdin-Francfort-Marigo~\cite{NumericalExperBourdi2000,NumericalImpleBourdi2007,TheVariationalBourdi2007,TheVariationalBourdi2008}, the total potential energy of the system with regularized crack description is expressed as
\begin{equation}
	\begin{aligned}
		\mathcal{E}_\epsilon(\bu, \alpha)
		\coloneqq \int_{\Omega_0}
		g(\alpha(\bx)) \Psi(\boldsymbol{\varepsilon}( & \bu(\bx)))
		\, \dd{\bx}
		+ \sum_{l=1}^L \int_{\Omega_l} \Psi(\boldsymbol{\varepsilon}(\bu(\bx))) \, \dd{\bx}                                                                                                                 \\
		                                              & + \frac{G_{\mathrm{c}}}{c_w} \int_{\Omega}  \left[ \frac{w(\alpha(\bx))}{\epsilon} + \epsilon {\vert \nabla \alpha(\bx) \vert}^2 \right]\, \dd{\bx}
		+ \int_{\Gamma_\mathrm{I}} \mathcal{G}(\jump{\bu(\bx)}) \, \dd{s},
	\end{aligned}
  \label{eq:potential_energy}
\end{equation}
where $\mathcal{G}$ represents the elastic energy density as a function of displacement jump at a point on an interface, defined as
\begin{align}
	\mathcal{G}(\jump{\bu(\bx)}) \coloneqq \frac{1}{2} k_\mathrm{I} \jump{\bu(\bx)} \cdot \jump{\bu(\bx)}.
\end{align}
The constitutive law for the interfaces presented above assumes isotropic elastic behavior. However, other complexities, such as anisotropic elasticity or cohesive behavior, can be incorporated into the model as long as the total potential energy functional~\eqref{eq:potential_energy} is twice differentiable.
Additionally, $G_\mathrm{c}$ denotes the fracture toughness of the matrix material, while $\epsilon > 0$ represents the numerical length parameter that characterizes the thickness of the damage zone around the crack.
The bulk fracture energy density is characterized by the crack geometric function $w(\alpha)$, which satisfies the conditions
\begin{align}
	w(0) = 0, w(1) = 1, w^\prime(\alpha) \geq 0\ \forall \alpha \in [0, 1].
\end{align}
Moreover, a normalization constant $c_w$ is defined as $c_w = 4 \int_0^1 \sqrt{w(t)}\,\dd{t}$ whose value is dependent on the choice of $w(\alpha)$.
The degradation function $g(\alpha)$ is a continuous monotonic function that fulfills the properties
\begin{align}
	g(0) = 1, g(1) = 0, g^\prime(\alpha) < 0 \ \forall \alpha \in [0, 1).
\end{align}
By having this specific choice for $w$ and $g$, the potential energy functional with the diffused crack $\Gamma$-converges~\cite{ApproximationOBraide1998,AnApproximatioChambo2004} to the potential energy functional with discrete crack for vanishing numerical length parameter $\epsilon$~\cite{NumericalExperBourdi2000,TheVariationalBourdi2008}.

\begin{figure}[ht]
	\centering
	\begin{subfigure}{0.49\textwidth}
		\begin{center}
			\includegraphics[width=\textwidth]{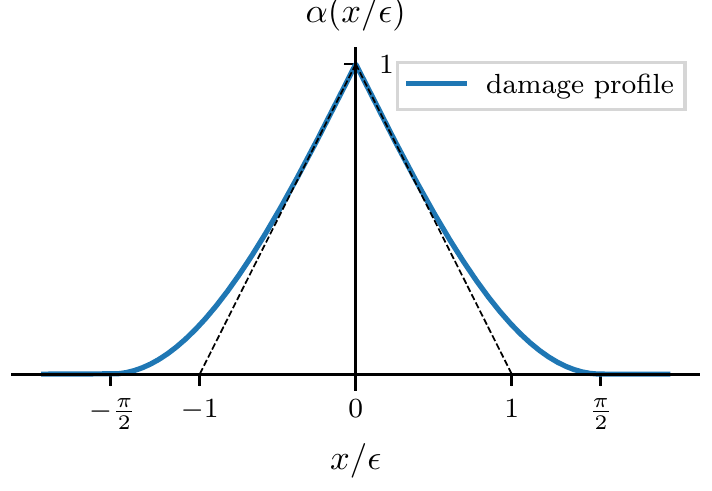}
		\end{center}
	\end{subfigure}
	\begin{subfigure}{0.49\textwidth}
	\def\svgwidth{1\columnwidth}
\begingroup%
  \makeatletter%
  \providecommand\color[2][]{%
    \errmessage{(Inkscape) Color is used for the text in Inkscape, but the package 'color.sty' is not loaded}%
    \renewcommand\color[2][]{}%
  }%
  \providecommand\transparent[1]{%
    \errmessage{(Inkscape) Transparency is used (non-zero) for the text in Inkscape, but the package 'transparent.sty' is not loaded}%
    \renewcommand\transparent[1]{}%
  }%
  \providecommand\rotatebox[2]{#2}%
  \newcommand*\fsize{\dimexpr\f@size pt\relax}%
  \newcommand*\lineheight[1]{\fontsize{\fsize}{#1\fsize}\selectfont}%
  \ifx\svgwidth\undefined%
    \setlength{\unitlength}{319.96796135bp}%
    \ifx\svgscale\undefined%
      \relax%
    \else%
      \setlength{\unitlength}{\unitlength * \real{\svgscale}}%
    \fi%
  \else%
    \setlength{\unitlength}{\svgwidth}%
  \fi%
  \global\let\svgwidth\undefined%
  \global\let\svgscale\undefined%
  \makeatother%
  \begin{picture}(1,0.62043951)%
    \lineheight{1}%
    \setlength\tabcolsep{0pt}%
    \put(0,0){\includegraphics[width=\unitlength,page=1]{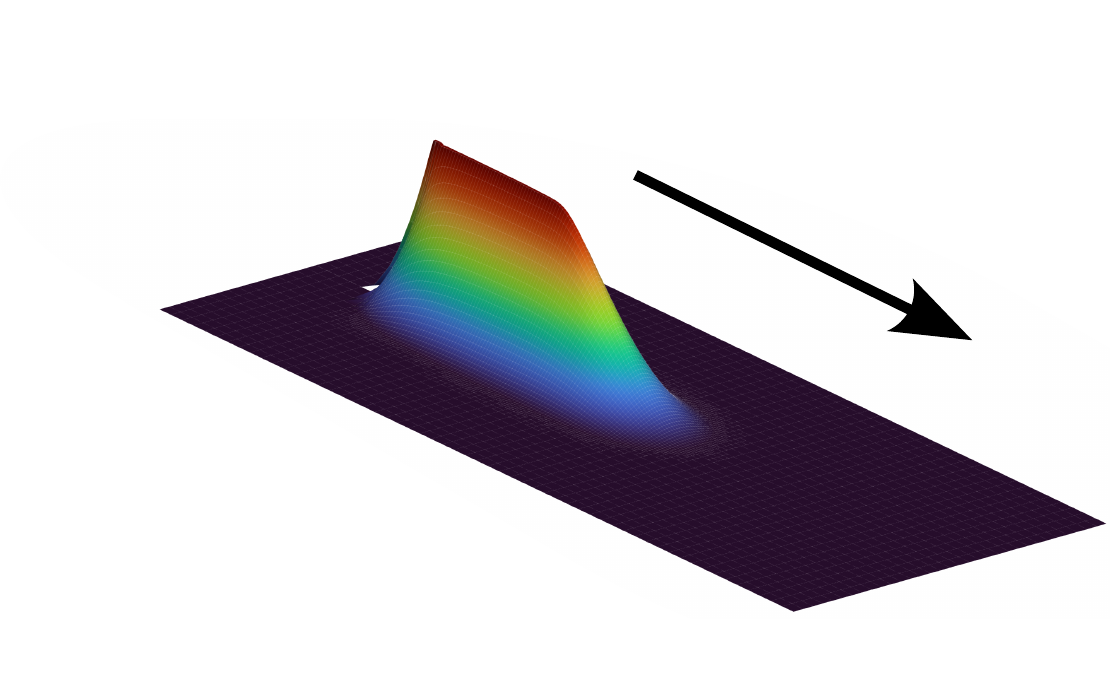}}%
    \put(0.74372547,0.43937502){\makebox(0,0)[lt]{\lineheight{1.25}\smash{\begin{tabular}[t]{l}\begin{tabular}{c}crack \\ propagation \\ direction \\\\\end{tabular}\end{tabular}}}}%
    \put(0,0){\includegraphics[width=\unitlength,page=2]{crack_profiles.pdf}}%
    \put(0.28307649,0.17836295){\makebox(0,0)[lt]{\lineheight{1.25}\smash{\begin{tabular}[t]{l}$x$\end{tabular}}}}%
    \put(0.28264445,0.06904331){\makebox(0,0)[lt]{\lineheight{1.25}\smash{\begin{tabular}[t]{l}$y$\end{tabular}}}}%
    \put(0,0){\includegraphics[width=\unitlength,page=3]{crack_profiles.pdf}}%
    \put(0.11873603,0.29752383){\makebox(0,0)[lt]{\lineheight{1.25}\smash{\begin{tabular}[t]{l}$\alpha$\end{tabular}}}}%
  \end{picture}%
\endgroup%

	\end{subfigure}
	\caption{
		Damage profile for the \texttt{PF-CZM} fracture model.
	}
	\label{fig:crack_profiles}
\end{figure}

Various models exist to define combinations of the crack geometric function $w(\alpha)$ and the degradation function $g(\alpha)$.
The most widely used models in the literature are \texttt{AT-2}~\cite{NumericalExperBourdi2000}, \texttt{AT-1}~\cite{MorphogenesisPropagationBourdin2014,GradientDamagePham2011}, and a more recently introduced, \texttt{PF-CZM}~\cite{AUnifiedPhaseWuJi2017,AGeometricallyWuJi2018,ALengthScaleWuJi2018} model.
While the \texttt{AT-2} and \texttt{AT-1} models are more suitable for modeling brittle fracture, the \texttt{PF-CZM} model combines the phase-field (PF) model with the cohesive zone model (CZM), enabling the use of a material traction-separation law.
The \texttt{PF-CZM} model allows definition of an independent characteristic length parameter $l_\mathrm{ch}$ (also called Irwin's internal length) which controls the fracture strength of the material with prescribed fracture toughness $G_\mathrm{c}$.
In this work, we adopt the \texttt{PF-CZM} model with linear softening traction separation law~\cite{ALengthScaleWuJi2018}, for which
the degradation function, crack geometric function and the normalization constant are given by
\begin{align}
	g(\alpha) \coloneqq \frac{{(1-\alpha)}^2}{{(1-\alpha)}^2 + \frac{4l_\mathrm{ch}}{\pi\epsilon} \alpha \left(1-\frac{\alpha}{2}\right)}, \qquad
  w(\alpha) \coloneqq 2\alpha-\alpha^2,
  \qquad
  \text{and}
  \qquad
c_w = \pi,
\end{align}
respectively.
The damage profile (Figure~\ref{fig:crack_profiles}) of the diffused crack is then given by
\begin{align}
	\alpha(x) = 1-\sin(\frac{x}{\epsilon})
\end{align}
with damage bandwidth of $\pi \epsilon$.

Following the principle of minimum potential energy of the mechanical system, we define the phase-field fracture problem as:
\begin{definition}[Potential Energy Minimization]
	\label{def:energy_minimization_pfm}
	Given a prescribed boundary displacement $\hat{\bu}^{(k)}$ on $\Gamma_\mathrm{D}$ at current time-step $k$ and phase-field $\alpha^{(k-1)}$ evaluated at the previous time step, find
	\begin{subequations}
		\begin{alignat}{3}
			\{ \bu^{(k)}, \alpha^{(k)} \} \in & \argmin_{\bu\in\mathcal{S}, \alpha\in\mathcal{A}} \  &  & \mathcal{E}_\epsilon(\bu, \alpha),                                                            \\
			                                  & \text{such that}\                                    &  & \bu = \hat{\bu}^{(k)} \quad \text{on } \Gamma_\mathrm{D}, & \qquad \text{(boundary condition)} \\
			                                  &                                                      &  & \alpha \geq \alpha^{(k-1)} \quad \text{in } \Omega.       & \qquad \text{(crack irreversibility)}
		\end{alignat}
	\end{subequations}
\end{definition}

The above minimization problem is non-convex in the unknown fields and incorporates bound constraints.
To write the strong form of the problem, we first have to write the Lagrangian of the constrained minimization problem, viz.,
\begin{equation}
	\begin{aligned}
		\mathcal{L}(\bu, \alpha, \lambda_1, \lambda_2)
		\coloneqq & \int_{\Omega_0} g(\alpha) \Psi(\boldsymbol{\varepsilon}) \, \dd{\bx}
		+ \sum_{l=1}^L \int_{\Omega_l} \Psi(\boldsymbol{\varepsilon}) \, \dd{\bx}
		+ \frac{G_{\mathrm{c}}}{c_w} \int_{\Omega}  \left[ \frac{w(\alpha)}{\epsilon} + \epsilon {\vert \nabla \alpha \vert}^2 \right]\, \dd{\bx}
		+ \int_{\Gamma_\mathrm{I}} \mathcal{G}(\jump{\bu}) \, \dd{s}                     \\
		          & - \int_\Omega \lambda_1 [\alpha - \alpha^{(k-1)}] \, \dd{\bx}
		+ \int_\Omega \lambda_2 [\alpha - 1] \, \dd{\bx},
	\end{aligned}
\end{equation}
where $\lambda_1: \Omega \rightarrow \mathbb{R}$ and $\lambda_2 : \Omega \rightarrow \mathbb{R}$ are the Lagrange multiplier fields corresponding to the inequality constraints on $\alpha$.
The first-order optimality conditions of the minimization problem are then given by
\begin{subequations}
	\begin{empheq}[left=\empheqlbrace]{align}
    \delta \mathcal{L}(\bu, \alpha, \lambda_1, \lambda_2, \delta\bu, \delta\alpha) = 0 \quad & \forall \delta\bu \in \mathcal{S},  \delta\alpha \in \mathcal{A}, \delta\bu = \boldsymbol{0} \text{ on } \Gamma_\mathrm{D}, \\
		\alpha^{(k-1)} \leq \alpha \leq 1 \quad                                                  & \text{in } \Omega ,                                                      \\
		\lambda_1 (\alpha - \alpha^{(k-1)}) = 0 \quad                                            & \text{in } \Omega ,                                                      \\
		\lambda_2 (\alpha - 1) = 0 \quad                                                         & \text{in } \Omega ,                                                      \\
		\lambda_1 \geq 0 \quad                                                                   & \text{in } \Omega ,                                                      \\
		\lambda_2 \geq 0 \quad                                                                   & \text{in } \Omega ,                                                      \\
		\bu = \hat{\bu}^{(k)} \quad                                                 & \text{on } \Gamma_\mathrm{D},
	\end{empheq}
\end{subequations}
where the variation of the Lagrangian with respect to $\bu$ and $\alpha$ is written as
\begin{equation}
	\begin{aligned}
		\delta \mathcal{L}(\bu, \alpha, & \lambda_1, \lambda_2, \delta\bu, \delta\alpha)                                                                                                                        \\
		=                               & \int_{\Omega_0}  g(\alpha) \bsigma : \delta \bepsilon \, \dd{\bx}
		+ \sum_{l=1}^L \int_{\Omega_l} \bsigma : \delta  \bepsilon \, \dd{\bx}
		+ \int_{\Omega_0} g^\prime(\alpha) \delta \alpha \Psi(\bepsilon) \, \dd{\bx}
		+ \int_{\Gamma_\mathrm{I}} \jump{\delta\bu} \cdot \btau(\jump{\bu}) \,\dd{s}                                                                                                                            \\
		                                & +  \frac{G_\mathrm{c}}{c_w} \int_{\Omega} \left[\frac{w^\prime(\alpha) \delta \alpha}{\epsilon} + 2 \epsilon \Grad \alpha \cdot \Grad \delta \alpha \right]\,\dd{\bx}
		- \int_\Omega \lambda_1 \delta \alpha \, \dd{\bx}
		+ \int_\Omega \lambda_2 \delta \alpha \, \dd{\bx}.
	\end{aligned}
\end{equation}
Here, $\btau(\jump{\bu}) = \nabla \mathcal{G}(\jump{\bu})$ is the interfacial traction.
To write the strong form of the problem, we use divergence theorem and assuming continuity of the traction along the interfaces, i.e.
\begin{align}
	\left.g(\alpha)\bsigma\cdot\bn \right\vert_{\Gamma_\mathrm{I}^+} =\left. \bsigma \cdot \bn\right\vert_{\Gamma_\mathrm{I}^-},
\end{align}
the above equation can be rewritten in the following form:
\begin{equation}
	\begin{aligned}
		\delta \mathcal{L}(\bu, \alpha, & \lambda_1, \lambda_2, \delta\bu, \delta\alpha)                                                                  \\
		=                               & -\int_{\Omega_0} \Div(g(\alpha) \bsigma) \cdot \delta \bu \, \dd{\bx}
		- \sum_{l=1}^L \int_{\Omega_l} \Div\bsigma \cdot \delta  \bu \, \dd{\bx}
		+ \int_{\Omega_0} g^\prime(\alpha) \delta \alpha \Psi(\bepsilon) \, \dd{\bx}                                                                      \\
		                                & + \int_{\Gamma_\mathrm{I}^-} \jump{\delta\bu} \cdot [\btau(\jump{\bu}) - \bsigma \cdot \bn] \,\dd{s}
		+  \frac{G_\mathrm{c}}{c_w} \int_{\Omega} \left[\frac{w^\prime(\alpha) }{\epsilon} - 2 \epsilon\Div \Grad \alpha   \right] \delta\alpha\,\dd{\bx} \\
		                                & - \int_\Omega \lambda_1 \delta \alpha \, \dd{\bx}
		+ \int_\Omega \lambda_2 \delta \alpha \, \dd{\bx} = 0.
	\end{aligned}
\end{equation}

Since the variations $\delta\bu$ and $\delta\alpha$ are arbitrary, the fracture problem can be expressed in strong form as:
\begin{definition}[Strong form]
	\label{def:pfm_strong_form}
	Given $\alpha^{(k-1)}$ in $\Omega$ and $\hat{\bu}^{(k)}$ on $\Gamma_\mathrm{D}$, find $\bu^{(k)} \equiv \bu : \overline{\Omega} \rightarrow \mathbb{R}^d$ and $\alpha^{(k)} \equiv \alpha : \Omega \rightarrow \mathbb{R}$ such that
	\begin{subequations}
		\begin{align}
			-  \Div (g(\alpha)\bsigma) = \boldsymbol{0} \quad                                                                                                                            & \text{in } \Omega_0 ,                \\
      - \Div \bsigma = \boldsymbol{0} \quad                                                                                                                                                     & \text{in } \Omega_l, l = 1,\dots,L , \\
			g^\prime(\alpha) \Psi(\bepsilon) + \frac{G_\mathrm{c}}{c_w} \left[ \frac{w^\prime(\alpha)}{\epsilon} - 2\epsilon \Div \Grad \alpha \right] - \lambda_1 + \lambda_2 = 0 \quad & \text{in } \Omega_0 ,                \\
			\frac{G_\mathrm{c}}{c_w} \left[ \frac{w^\prime(\alpha)}{\epsilon} - 2\epsilon \Div \Grad \alpha \right] - \lambda_1 + \lambda_2 = 0 \quad                                    & \text{in } \Omega_l, l=1,\dots,L ,   \\
			g(\alpha)\bsigma \cdot \bn = \btau(\jump{\bu}) \quad                                                                                                                         & \text{on } \Gamma_\mathrm{I}^+ ,     \\
			\bsigma \cdot \bn = \btau(\jump{\bu}) \quad                                                                                                                                  & \text{on } \Gamma_\mathrm{I}^- ,     \\
			\alpha^{(k-1)} \leq \alpha \leq 1 \quad                                                                                                                                      & \text{in } \Omega ,                  \\
			\lambda_1 [\alpha - \alpha^{(k-1)}] = 0 \quad                                                                                                                                & \text{in } \Omega ,                  \\
			\lambda_2 [\alpha - 1] = 0 \quad                                                                                                                                             & \text{in } \Omega ,                  \\
			\lambda_1 \geq 0 \quad                                                                                                                                                       & \text{in } \Omega ,                  \\
			\lambda_2 \geq 0 \quad                                                                                                                                                       & \text{in } \Omega ,                  \\
			\alpha = 0 \quad                                                                                                                                                             & \text{on } \partial\Omega ,          \\
			\bsigma \cdot \bn = \boldsymbol{0} \quad                                                                                                                                     & \text{on } \Gamma_\mathrm{N},        \\
			\bu = \hat{\bu}^{(k)} \quad                                                                                                                                                  & \text{on } \Gamma_\mathrm{D} .
		\end{align}
	\end{subequations}
\end{definition}

In the next sections, we present the spatial discretization scheme and a solution algorithm to solve the phase-field fracture problem.

\section{Spatial discretization}
\label{sec:spatial_discretization}

In this study, we utilize the extended finite element method (XFEM)~\cite{TheExtendedGeFries2010,OnImplementationXFEMCarraro2015} to simulate displacement discontinuities that occur at interfaces. XFEM is a versatile numerical method capable of handling a wide range of interface problems. To ensure comprehensiveness, we present a concise summary of its fundamental formulation in the following.

To partition the spatial domain $\Omega$ into subdomains $\Omega_l,\ l=1,\dots,L$, we utilize smooth level-set functions $\psi_l : \Omega \rightarrow \mathbb{R}$ that satisfy the condition $\psi_l = 0$ on the corresponding interfaces.
Given a discretized domain $\Omega_h$ using finite elements, let $\omega_{h,l}^\mathrm{cut} \subset \Omega_h$ be a subdomain covered by all the cells cut by an interface represented by the level-set function $\psi_l$, and $\omega_{h}^\mathrm{uncut} \coloneqq \Omega_h \setminus \cup_{l=1}^{L} \omega_{h,l}^\mathrm{cut}$ be the part of the spatial domain covered by all the non-intersected cells.
At time-step $k$, we express the discretized displacement field as
\begin{equation}
	\begin{aligned}
		\bu^{(k)}(\bx) \approx \bu_h(\bx; \ba^{(k)}, \bb^{(k)}) \coloneqq
		\begin{cases}
			\sum_{i \in I_a}  a_i^{(k)} \bphi_i(\bx)                                                          & \text{for } \bx \in \omega_h^\mathrm{uncut},   \\
			\sum_{i \in I_a}  a_i^{(k)} \bphi_i(\bx) + \sum_{i \in I_b} b_i^{(k)} \xi_{l,i}(\bx) \bphi_i(\bx) & \text{for } \bx \in \omega_{h,l}^\mathrm{cut},
		\end{cases}
	\end{aligned}
	\label{eq:discretized_displacement_field}
\end{equation}
where $\bphi_i$ are the conventional finite element shape functions.
The field vectors $\ba^{(k)} \in \mathbb{R}^{|I_a|}$ and $\bb^{(k)} \in \mathbb{R}^{|I_b|}$ represent conventional and enriched nodal displacements, respectively, where $I_a$ and $I_b$ are the corresponding index sets.
To model the discontinuities within the cut cells, a usual XFEM approach uses an enrichment function of the form
\begin{equation}
	\xi_{l,i}(\bx) \coloneqq \sgn(\psi_l (\bx)) - \sgn(\psi_l (\boldsymbol{X}_i)),
	\label{eq:discontinuous_enrichment_function}
\end{equation}
where $\boldsymbol{X}_i$ is the position of the finite element node corresponding to the $i^\mathrm{th}$ degree of freedom.
The displacement jump at a point $\bx \in \Gamma_\mathrm{I}$ is given by
\begin{align}
  \jump{\bu^{(k)}(\bx)} \approx  \jump{\bu_h(\bx; \ba^{(k)}, \bb^{(k)})} = 2\sum_{i \in I_b} b_i^{(k)} \bphi_i(\bx).
\end{align}
For more details on the implementation of XFEM for interface problems, please refer to~\cite{OnImplementationXFEMCarraro2015}.

Since the phase-field $\alpha$ does not admit discontinuities at the interfaces, we discretize it using the standard finite element approach.
At time-step $k$, we write the discretized phase field as
\begin{align}
	\alpha^{(k)}(\bx) \approx \alpha_h(\bx; \boldsymbol{c}^{(k)})
	\coloneqq \sum_{i \in I_c} c_i^{(k)} \varphi_i(\bx) \quad \forall \bx \in \Omega_h,
\end{align}
where $\varphi_i$ are the scalar finite element shape functions.
The field vector $\boldsymbol{c}^{(k)} \in \mathbb{R}^{\vert I_c \vert}$ represents nodal values for phase field, where $I_c$ denotes the corresponding index set.

Throughout this paper, all computational investigations were conducted using the open-source finite element library, \texttt{deal.II}~\cite{TheDealIiFinArndt2021}, in which bilinear finite elements were utilized for both the displacement and phase fields in the finite-element problem implementation.

\section{Simulation algorithm}\label{sec:simulation_algorithm}

Since the potential energy minimization problem (Definition~\ref{def:energy_minimization_pfm}) is twice-differentiable with respect to the unknown fields $\bu$ and $\alpha$, it can be solved using a continuous, gradient-based optimization solver.
For the numerical studies presented in this paper, we used \texttt{IPOPT} ~\cite{OnTheImplemenWachte2006} (version 3.14.9),
which is an open-source software for solving large scale nonlinear optimization problems.
To get a discrete-time solution that approximates a quasi-static process, the optimizer is given a warm start with the solution either from the previous time step or that obtained by a transfer from old mesh to new mesh after a mesh refinement or coarsening cycle.
The minimization problem is solved until convergence, with monotonic decrease of the slack variable (\texttt{mu\_strategy}) with initial value \num{1e-9} (\texttt{mu\_init}).
All the other parameters of \texttt{IPOPT} are kept as default values.


In order to improve the efficiency and reliability of the solution algorithm for the phase-field fracture problem, which deals with periodic heterogeneous microstructures under surfing boundary conditions, we propose three additional techniques:
\emph{adaptive mesh refinement}, \emph{presolving} and \emph{backtracking}.
There steps are explained in the following subsections, while the complete solution procedure is presented as pseudo-code in Algorithm~\ref{alg:phase_field_simulation_algorithm}.

\subsection{Adaptive mesh refinement}

\begin{figure}[ht]
	\centering
	\def\svgwidth{1\columnwidth}
\begingroup%
  \makeatletter%
  \providecommand\color[2][]{%
    \errmessage{(Inkscape) Color is used for the text in Inkscape, but the package 'color.sty' is not loaded}%
    \renewcommand\color[2][]{}%
  }%
  \providecommand\transparent[1]{%
    \errmessage{(Inkscape) Transparency is used (non-zero) for the text in Inkscape, but the package 'transparent.sty' is not loaded}%
    \renewcommand\transparent[1]{}%
  }%
  \providecommand\rotatebox[2]{#2}%
  \newcommand*\fsize{\dimexpr\f@size pt\relax}%
  \newcommand*\lineheight[1]{\fontsize{\fsize}{#1\fsize}\selectfont}%
  \ifx\svgwidth\undefined%
    \setlength{\unitlength}{498.1594344bp}%
    \ifx\svgscale\undefined%
      \relax%
    \else%
      \setlength{\unitlength}{\unitlength * \real{\svgscale}}%
    \fi%
  \else%
    \setlength{\unitlength}{\svgwidth}%
  \fi%
  \global\let\svgwidth\undefined%
  \global\let\svgscale\undefined%
  \makeatother%
  \begin{picture}(1,0.35989663)%
    \lineheight{1}%
    \setlength\tabcolsep{0pt}%
    \put(0,0){\includegraphics[width=\unitlength,page=1]{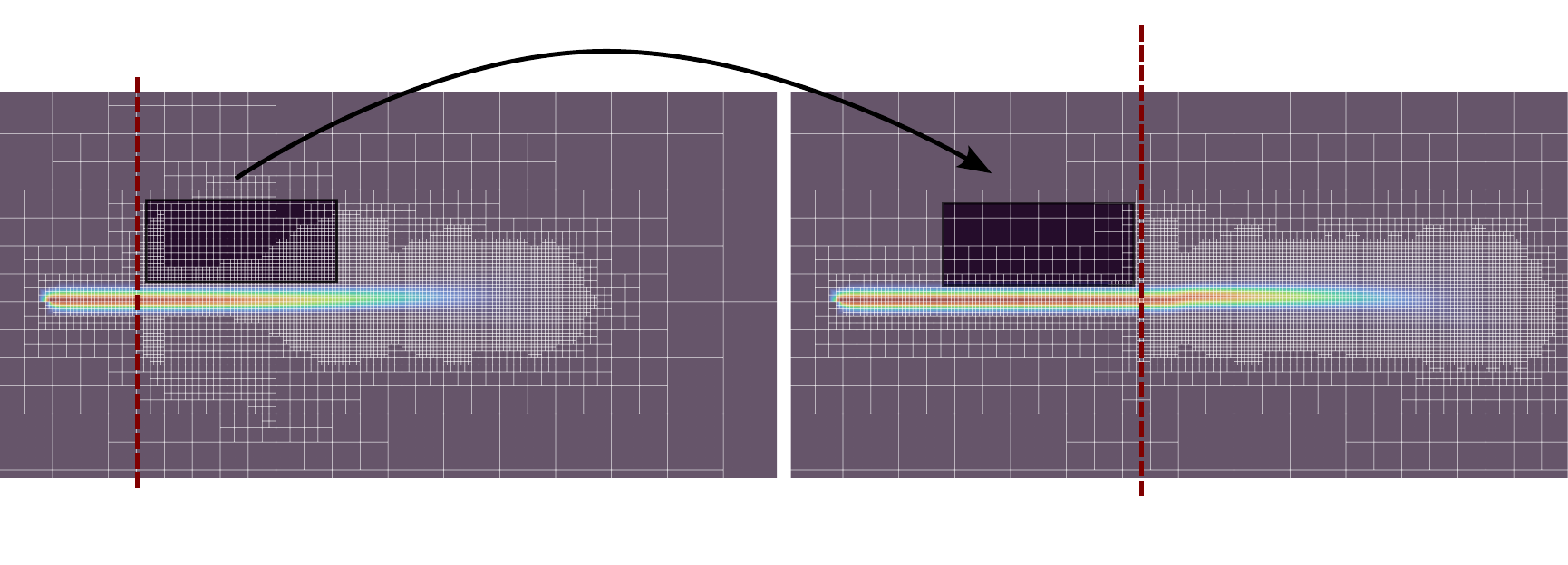}}%
    \put(0.34340137,0.34044501){\color[rgb]{0,0,0}\makebox(0,0)[lt]{\lineheight{1.25}\smash{\begin{tabular}[t]{l}Coarsening\end{tabular}}}}%
    \put(0,0){\includegraphics[width=\unitlength,page=2]{adaptive_mesh.pdf}}%
    \put(0.56918298,0.00100864){\color[rgb]{0,0,0}\makebox(0,0)[lt]{\lineheight{1.25}\smash{\begin{tabular}[t]{l}Refinement\end{tabular}}}}%
    \put(0,0){\includegraphics[width=\unitlength,page=3]{adaptive_mesh.pdf}}%
    \put(0.5265734,0.34541754){\color[rgb]{0,0,0}\makebox(0,0)[lt]{\lineheight{1.25}\smash{\begin{tabular}[t]{l}Coarsening zone\end{tabular}}}}%
    \put(0.7479936,0.34541754){\color[rgb]{0,0,0}\makebox(0,0)[lt]{\lineheight{1.25}\smash{\begin{tabular}[t]{l}Refinement zone\end{tabular}}}}%
    \put(0.15771731,0.02450631){\color[rgb]{0,0,0}\makebox(0,0)[lt]{\lineheight{1.25}\smash{\begin{tabular}[t]{l}Crack tip\end{tabular}}}}%
    \put(0,0){\includegraphics[width=\unitlength,page=4]{adaptive_mesh.pdf}}%
  \end{picture}%
\endgroup%

	\caption{
		Illustration of adaptive mesh refinement and coarsening strategy based on damage thresholds and horizontal position of the crack tip.
	}
	\label{fig:adaptive_mesh}
\end{figure}

In a structural optimization framework, it is necessary to efficiently solve multiple state problems for various design points. To accomplish this, we employ a straightforward adaptive mesh refinement and coarsening strategy specifically to resolve the phase-field cracks. Additionally, under the surfing boundary conditions, the crack typically propagates in the direction of the surfing velocity (from left to right in the present case), thus we can incorporate mesh coarsening behind the crack tip to further minimize the simulation costs.
The strategy can be divided into the following components:
\begin{itemize}
	\item \emph{Finding crack tip position.} In this step, we traverse through all the cells in the finite element mesh and locate the maximum value of coordinate in the direction of the surfing velocity that has a phase-field value of $\alpha=1$.
	      This gives an approximate location of the crack tip.
	\item \emph{Refinement.} Mark cells for refinement ahead of the crack tip, wherever $\alpha > 0$.
	      Since we use the \texttt{PF-CZM} fracture model with finite damage bandwidth, the need for mesh refinement is limited to a small region around the crack.
	      To increase the number of refined cells around the crack tip, we also refine the neighboring cells and also their neighbors.
	\item \emph{Coarsening.} Mark cells for coarsening behind the crack tip with some margin, with a predefined phase-field threshold (for instance, $\alpha^\mathrm{th}=0.8$ as used in this work).
	      This step is performed only once in every time step to eliminate potential oscillations between refinement and coarsening.
\end{itemize}
These steps are also illustrated in Figure~\ref{fig:adaptive_mesh}.

\subsection{Presolving}
Prior to solving the coupled optimization problem (Definition~\ref{def:energy_minimization_pfm}) for unknown displacement and phase field,
the optimizer is initialized with the solution from the previous time step, and the constraints on the displacement degrees of freedom are updated according to the current time step.
It has been observed that if there is a large change in the boundary displacements between two consecutive time steps, damage may initiate in the finite elements located near the boundary.
Additionally, if a mesh refinement or coarsening step has been performed before solving the optimization problem, the solution provided as an initial guess may not be sufficiently accurate. This issue is particularly prevalent when enriched finite elements exhibiting strong displacement discontinuities are added or removed from the overall finite element mesh.

For these reasons, before solving the fully coupled optimization problem, we perform a presolve for the displacement field while keeping the phase field fixed.
Since the energy functional is strictly convex in the displacement variable, the uniqueness of the solution is ensured.
Following the update of the displacements for the current boundary conditions or current mesh, the energy minimization problem is then solved for both unknowns simultaneously.

\subsection{Backtracking}
\begin{figure}[htbp]
	\begin{subfigure}{0.49\textwidth}
		\centering
		\begin{tikzpicture}[text=black, draw=black]
			\begin{axis}[
					width=\textwidth,
					grid=both,
					xlabel={$t$},
					ylabel={J-integral},
					xlabel near ticks,
					ylabel near ticks,
					xmin=0, xmax=1.,
				]
				\addplot[mark=o, draw=blue, mark size=1pt] table[header=false, x index=1, y index=3] {data/backtracking/stats_without_backtracking.dat};
				\addplot[mark=o, draw=red, mark size=1pt] table[header=false, x index=1, y index=3] {data/backtracking/stats_with_backtracking.dat};
				\node at (axis cs:0.53,2.85) {(a)};
				\node at (axis cs:0.5,0.95) {(b)};
				\node at (axis cs:0.2,1.1) {(c)};
			\end{axis}
		\end{tikzpicture}
	\end{subfigure}
	\begin{subfigure}{0.49\textwidth}
		\centering
		\begin{tikzpicture}[text=black, draw=black]
			\begin{axis}[
					width=\textwidth,
					grid=both,
					xlabel={$t$},
					ylabel={Crack length},
					xlabel near ticks,
					ylabel near ticks,
					xmin=0, xmax=1.,
					legend style={at={(0.03,1)},anchor=north west}
				]
				\addplot[mark=o, draw=blue, mark size=1pt] table[header=false, x index=1, y index=9] {data/backtracking/stats_without_backtracking.dat};
				\addlegendentry{without backtracking};
				\addplot[mark=o, draw=red, mark size=1pt] table[header=false, x index=1, y index=9] {data/backtracking/stats_with_backtracking.dat};
				\addlegendentry{with backtracking};
				\node at (axis cs:0.48,8) {(a)};
				\node at (axis cs:0.48,37) {(b)};
				\node at (axis cs:0.17,19) {(c)};
			\end{axis}
		\end{tikzpicture}
	\end{subfigure}
	\caption{
		Illustration of the concept of backtracking step in phase-field fracture simulation.
		The left plot displays the J-integral curves, a measure of fracture resistance, while the right plot shows the corresponding macroscopic crack lengths for simulations both with and without backtracking step.
		The backtracking step is demonstrated by reversing time in the event that the crack advance exceeds a predefined range, as shown by the progression from (a) to (b).
		After attaining the crack advance within an acceptable range (a)-(c), the time is again set to move forward.
		The crack states corresponding to points (a), (b) and (c) are depicted in Figure~\ref{fig:crack_backtracking}.
	}
  \label{fig:backtracking}
	\bigskip
	\begin{center}
		\begin{subfigure}{\textwidth}
			\centering
			\includegraphics[width=0.49\textwidth]{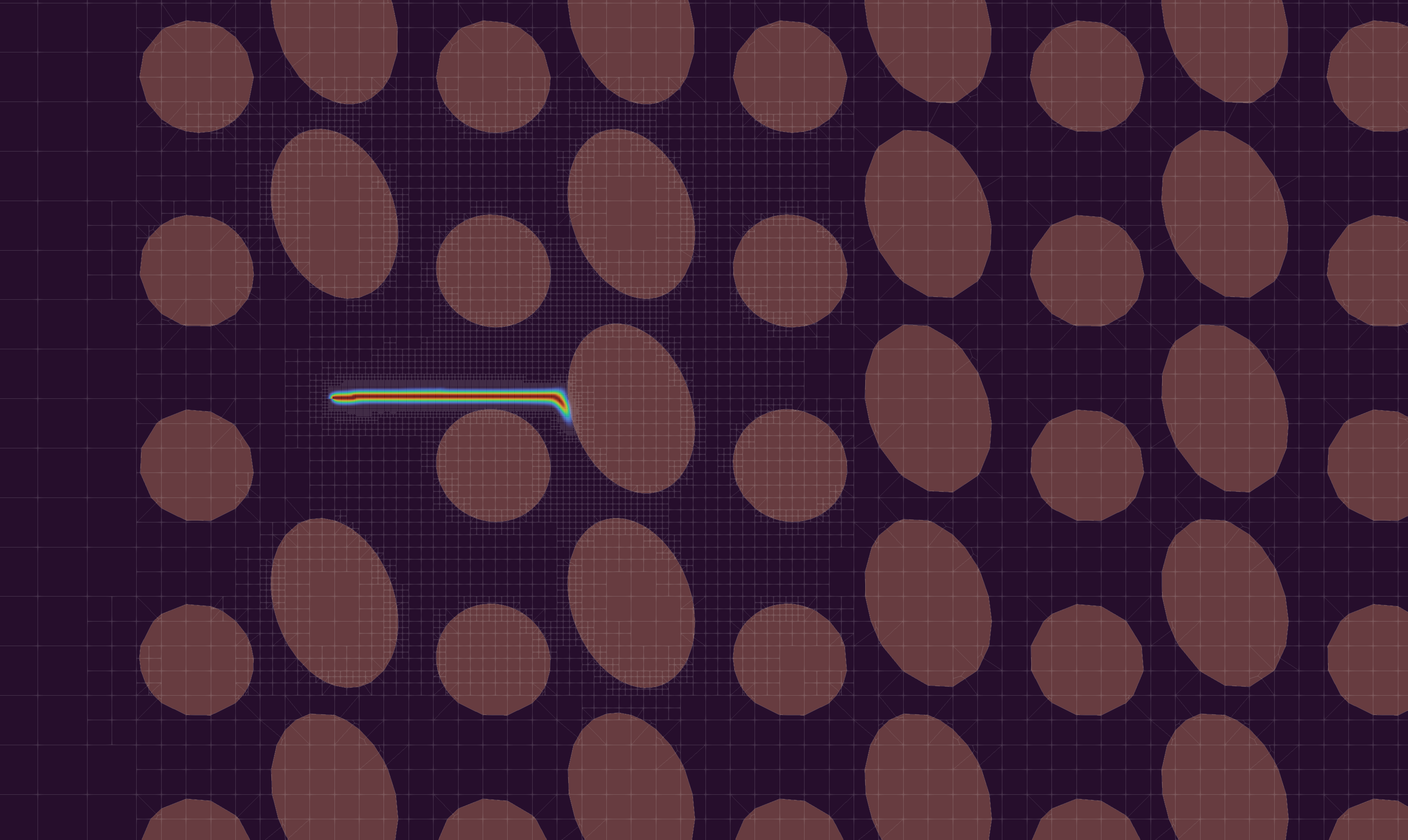}
			\caption{Crack at previous time step}
		\end{subfigure}
		\medskip
		\begin{subfigure}{0.49\textwidth}
			\includegraphics[width=\textwidth]{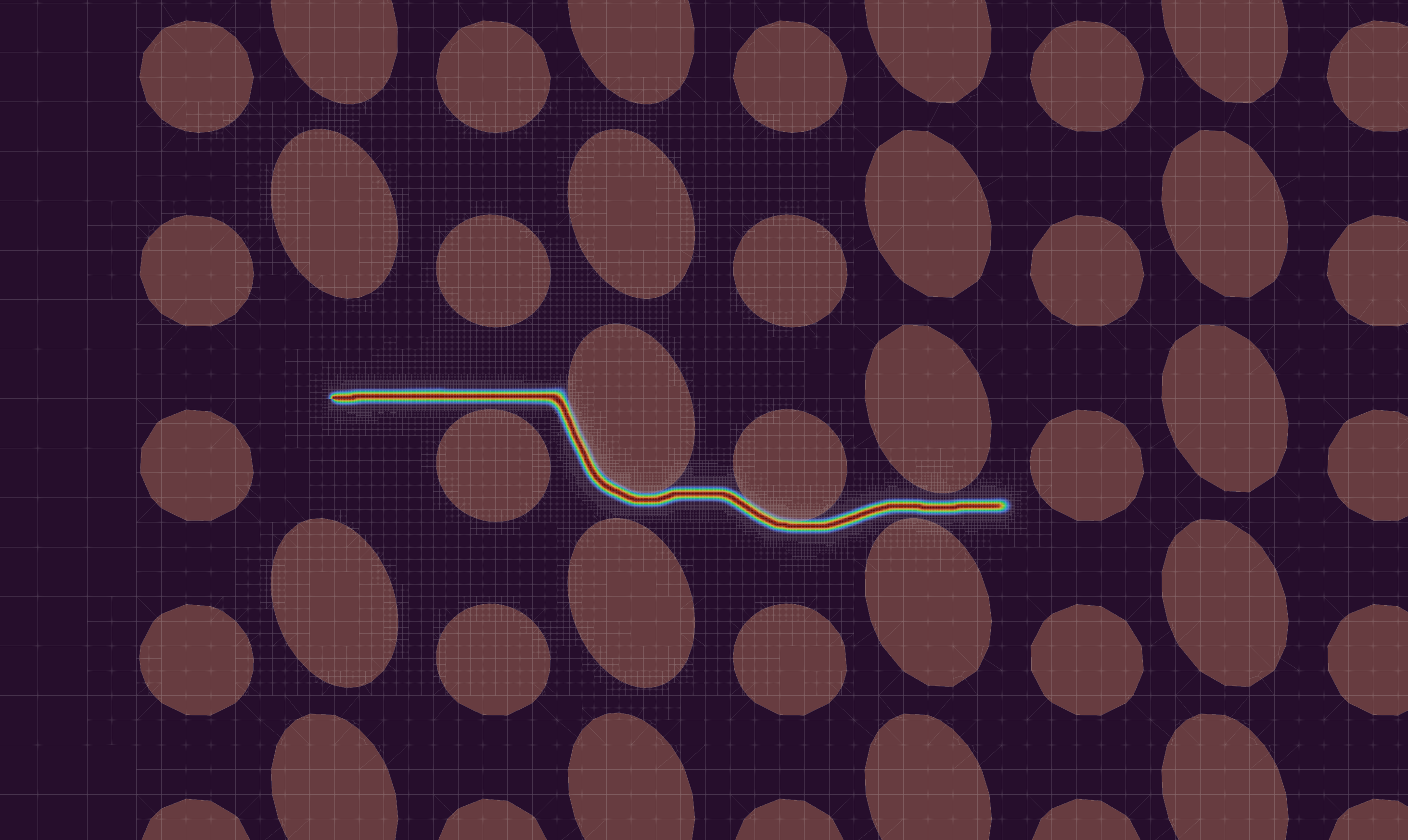}
			\caption{Crack without backtracking}
		\end{subfigure}
		\begin{subfigure}{0.49\textwidth}
			\includegraphics[width=\textwidth]{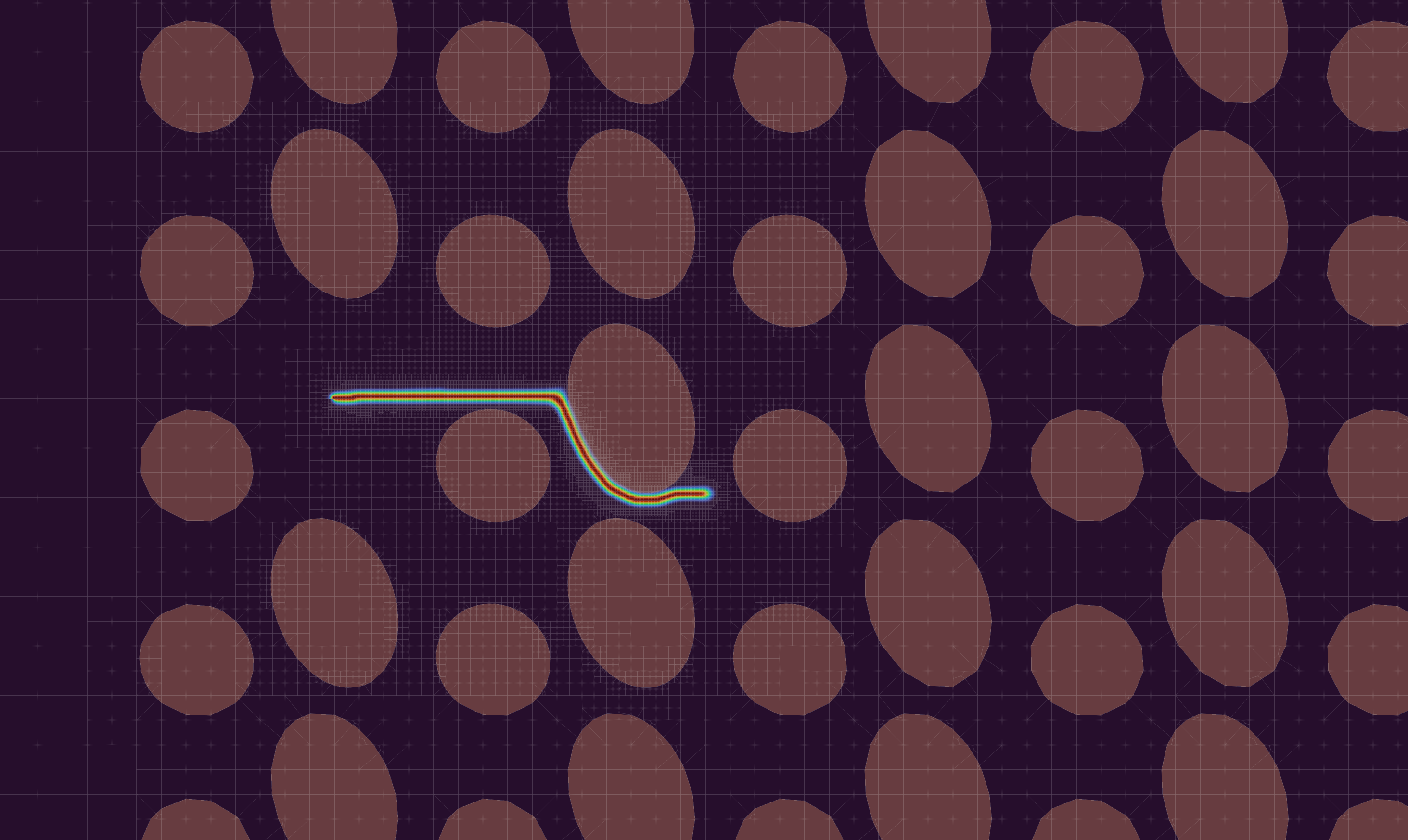}
			\caption{Crack with backtracking}
		\end{subfigure}
	\end{center}
	\caption{Illustration of the effect of using backtracking step in phase-field fracture simulation.}
	\label{fig:crack_backtracking}
\end{figure}

If the crack progresses too far in a single time step, it may extend beyond a representative volume element (RVE) in a periodic microstructure, potentially leading to the omission of important crack resistance behavior within that RVE.
To address this issue, it may be necessary to implement a backtracking procedure in order to ensure that the crack advances at a reasonable rate. This involves resetting the time to a previous point in the simulation and solving the minimization problem (Definition~\ref{def:energy_minimization_pfm}) again using the current solution as an initial guess. The crack tip is expected to move in the reverse direction, and the procedure is repeated until the advance of the crack tip between time steps is deemed acceptable.
The concept is illustrated in Figure~\ref{fig:backtracking} and Figure~\ref{fig:crack_backtracking}.

It is worth noting that a similar backtracking approach has been used in the past for the purpose of avoiding local minima in the total potential energy functional~\cite{TheVariationalBourdi2008} and for gaining numerical stability of topology optimization algorithm for brittle fracture~\cite{TopologyOptimiDesai2022}. However, the goal in the present context is not to escape from local minima or improve numerical stability, but rather to limit the advance of the crack tip between time steps in order to ensure that the fracture resistance behavior of the RVE is properly captured.

Since the crack tip travels in the reverse direction during a backtracking loop, no local mesh refinement and coarsening step is performed at this stage.

\begin{algorithm2e}[!h]
	\SetAlgoLined
	\DontPrintSemicolon
	\caption{Phase-field fracture simulation with surfing boundary conditions}\label{alg:phase_field_simulation_algorithm}
	\KwIn{Time-step size: $\Delta t$, max allowed crack advance, target crack length}
	\KwResult{Solution: $\boldsymbol{a}^{(k)}$, $\boldsymbol{b}^{(k)}$, $\boldsymbol{c}^{(k)}$, $k=0,\dots,N_T$}

	$k \gets 0$ \;
	$t \gets 0$ \;
	$\boldsymbol{a}^{(k)}=\boldsymbol{0},\boldsymbol{b}^{(k)}=\boldsymbol{0},\boldsymbol{c}^{(k)}=\boldsymbol{0}, \boldsymbol{c}^{(k-1)}=\boldsymbol{0}$ \;
	$\tilde{\bx}^{(-1)} \gets$ initial crack tip position \;
	coarsen flag $\gets 1$ \;
	\While{$t \leq T$ \textbf{or} crack length < target crack length} {

	construct boundary constraints for time $t$ and hanging-node constraints \;
	solve for displacement field only \hspace{5cm}
	$ \boldsymbol{a}^{(k)}, \boldsymbol{b}^{(k)} \gets \argmin_{\boldsymbol{a}, \boldsymbol{b}} \mathcal{E}_\epsilon(\bu_h(\,\cdot\,,\boldsymbol{a}, \boldsymbol{b}), \alpha_h(\,\cdot\,, \boldsymbol{c}^{(k)})) $ \;
	solve for displacement and phase field with irreversibility constraints
	$\boldsymbol{a}^{(k)}, \boldsymbol{b}^{(k)}, \boldsymbol{c}^{(k)} \gets \argmin_{\boldsymbol{a}, \boldsymbol{b}, \boldsymbol{c}} \mathcal{E}_\epsilon(\bu_h(\,\cdot\,,\boldsymbol{a}, \boldsymbol{b}), \alpha_h(\,\cdot\,, \boldsymbol{c}))$
	s.t.\ $c_i^{(k)} \geq c_i^{(k-1)}\ \forall i \in I_d$ \;
	$\tilde{\bx}^k \gets$  find crack tip position \;
	\While{crack advance ($\tilde{x}_1^{(k)} - \tilde{x}_1^{(k-1)}$) $>$ max allowed crack advance}{
	$t \gets t - \Delta t$ \tcp*{Decrement time}
	construct boundary constraints for time $t$ and hanging-node constraints \;
	solve for displacement and phase field with irreversibility constraints
	$\boldsymbol{a}^{(k)}, \boldsymbol{b}^{(k)}, \boldsymbol{c}^{(k)} \gets \argmin_{\boldsymbol{a}, \boldsymbol{b}, \boldsymbol{c}} \mathcal{E}_\epsilon(\bu_h(\,\cdot\,,\boldsymbol{a}, \boldsymbol{b}), \alpha_h(\,\cdot\,, \boldsymbol{c}))$ s.t.\ $c_i^{(k)} \geq c_i^{(k-1)}\ \forall i \in I_d $ \;
	$\tilde{\bx}^{(k)} \gets$  find crack tip position \;
	}

	mesh refined or coarsened $\gets$ false \;
	\If{coarsen flag \textbf{is} $1$}{
		mark cells for coarsening \;
		coarsen flag $\gets 0$ \;
	}

	\While{true}{
		mark cells for refinement \;
		\uIf{cells due for refinement or coarsening}{
			perform local mesh refinement and coarsening \;
			project $\boldsymbol{c}^{(k)}$ and $\boldsymbol{c}^{(k-1)}$ to new mesh \;
			mesh refined or coarsened $\gets$ true
		}\Else{
			\Break \;
		}
	}
	\If{mesh refined or coarsened}{
		\Continue \;
	}

	coarsen flag $\gets 1$ \;

	$t \gets t + \Delta t$ \tcp*{Increment time}
	$k \gets k + 1$ \tcp*{Increment time step count}
	}
\end{algorithm2e}

\section{Microstructure optimization}\label{sec:optimization}

\begin{figure}[ht]
	\begin{center}
		\begin{tikzpicture}[scale=0.8, text=black, draw=black]

	\def\h{-0.15}
	\def\hx{0.0}
	\def\hy{0.2}
	\draw[fill=gray, fill opacity=0.05, draw=none, thick] (-1, -3) rectangle (11, 3);

	\foreach \a in {-1.,0.5,...,11} {
			\foreach \b in {-3, -2,...,3} {
					\draw[fill=teal!50, rotate around={45:(\a, \b+\h)}] (\a, \b+\h) ellipse (0.25 and 0.4);
				}
		}
	\foreach \a in {-0.25,1.25,...,11} {
			\foreach \b in {-2.5, -1.5,...,3} {
					\draw[fill=brown!30, rotate around={20:(\a+\hx, \b+\h+\hy)}] (\a+\hx, \b+\h+\hy) ellipse (0.3 and 0.45);
				}
		}
	\draw[draw=none, fill=white, opacity=1] (-2, 3) rectangle (12, 4);
	\draw[draw=none, fill=white, opacity=1] (-2, -3) rectangle (12, -4.);
	\draw[draw=none, fill=white, opacity=1] (11, -3.1) rectangle (12, 3.1);
	\draw[draw=none, fill=white, opacity=1] (-2, -3.1) rectangle (-1, 3.1);
	\draw[draw=black, thick, fill=gray, fill opacity=0.2] (-2, -4) rectangle (12, 4);
	\draw[draw=black, draw opacity=0.5, thick, dashed] (-1, -3) rectangle (11, 3);

	\draw[fill=white, draw=none] (-2.03, 0) rectangle (0, 0.03);
	\draw (-2.005, 0) -- (0, 0) -- (0, 0.03) -- (-2.005, 0.03);
	\node[anchor=south west] at (-2.05, 0) {Initial crack};

	\draw[red, thick] (2.75, 0.55) rectangle (4.25, 1.55) node[anchor=south, red] {RVE};

	\dimline[label style={above, fill=none}, extension start length=0.9cm, extension end length=0.9cm] {(4.25, 3.2)}{(5.75, 3.2)}{$15$};
	\dimline[label style={above, fill=none}, extension start length=0.6cm, extension end length=0.6cm] {(2, 3.2)}{(3.5, 3.2)}{$15$};
	\dimline[label style={above, fill=none}, extension start length=0.2cm, extension end length=4.3cm] {(-2, 4.3)}{(0, 4.3)}{$20$};
	\dimline[label style={above, fill=none}, extension start length=0.cm, extension end length=0cm] {(-2, 2)}{(-1, 2)}{$10$};
	\dimline[label style={above, fill=none}, extension start length=0.cm, extension end length=0cm] {(11, 2)}{(12, 2)}{$10$};
	\dimline[label style={above, fill=none}, extension start length=0.cm, extension end length=0cm] {(10, -4)}{(10, -3)}{$10$};
	\dimline[label style={above, fill=none}, extension start length=0.cm, extension end length=0cm] {(10, 3)}{(10, 4)}{$10$};
	\dimline[label style={above, fill=none}, extension start length=0.cm, extension end length=0.2cm] {(0, 4.3)}{(12, 4.3)}{$120$};

\end{tikzpicture}
	\end{center}
	\caption{
		Geometrical setup of heterogeneous microstructure for the design optimization problem.
		The two sets of inclusions are arranged in a periodic fashion such that orthotropic material behavior is observed at the macroscale.
		A fixed initial crack is provided on the left edge of the geometry while the whole inclusion setup is allowed to move vertically relative to the initial crack.
		The representative volume element (RVE) is illustrated by the red box.
	}
	\label{fig:geometrical_setup_optimization}
  \bigskip
	\centering
	\footnotesize{
	\def\svgwidth{1\columnwidth}
	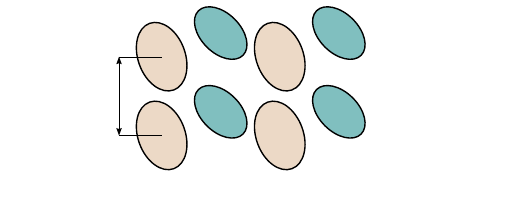

	}
  \caption{
    Design parameterization for the optimization problem with $d=9$ design parameters: $x_1,\dots,x_9$.
  }
	\label{fig:paper3_design_parameters}
\end{figure}

The objective of this research is to determine the optimal spatial arrangement of stiff inclusions embedded in a matrix material to enhance resistance to fracture caused by the propagation of matrix cracks. To accomplish this, we consider a microstructure as depicted in Figure~\ref{fig:geometrical_setup_optimization}, which comprises of two sets of elliptical inclusions arranged in a periodic manner within the matrix. The inclusions are allowed to vary in size, rotate and shift positions, as long as there exists a minimum clearance between the inclusions and maintain a periodic pattern.

\begin{table}[ht]
  \caption{Parameters for phase-field fracture problem.}
  \begin{subtable}{0.5\textwidth}
	\caption{Material parameters}
	\label{tab:material_parameters}
	\begin{center}
		\begin{tabular}[c]{l|c|c|c}
			\toprule
			Property        & Matrix & Inclusion & Interface                 \\
			\midrule
			\midrule
		  Elastic modulus & $1$    & $5$       & $k_\mathrm{I} =\num{100}$ \\
			Poisson's ratio   & $0.3$  & $0.3$     & $-$                       \\
      Characteristic length $l_\mathrm{ch}$                                                      & \num{1} &  $-$ & $-$                           \\
      Fracture toughness $G_\mathrm{c}$                                                   & $1$  & $-$ & $-$                               \\
			\bottomrule
		\end{tabular}
	\end{center}
  \end{subtable}
  \begin{subtable}{0.55\textwidth}
	\caption{Algorithmic parameters}
	\label{tab:algorithmic_parameters}
	\begin{center}
		\begin{tabular}[c]{l|c}
			\toprule
			Parameter                                                                                  & Value                              \\
			\midrule
			\midrule
      Mode-\rom{1} stress intensity factor $K_\mathrm{I}$ (Eq.~\eqref{eq:surfing_boundary_conditions}) & \num{1.1}                          \\
			Surfing velocity $v$ (Eq.~\eqref{eq:polar_coordinates_surfing})                                    & \num{50}                           \\
      Time-step size $\Delta t$                                                                             & \num{0.01}                         \\
			Mesh size to numerical length ratio $h/\epsilon$                                           & $0.4$                              \\
			Damage threshold for local mesh refinement                                                 & \num{0.001}                        \\
			Damage threshold for local mesh coarsening                                                 & \num{0.8}                          \\
			Damage threshold for irreversibility                                                       & \num{0.5}                          \\
			\bottomrule
		\end{tabular}
	\end{center}
  \end{subtable}
  \bigskip
	\caption{Design parameters and bounds (cf.~Figure~\ref{fig:paper3_design_parameters})}.
	\label{tab:design_parameters}
	\begin{center}
		\begin{tabular}[c]{l|c|c|c}
			\toprule
			Parameter                             & Symbol     & Lower bound & Upper bound \\
			\midrule
			\midrule
			Horizontal alignment                  & $x_1$      & $0$         & $0.5$       \\
			Vertical alignment                    & $x_2$      & $0$         & $0.5$       \\
			Ellipse radius \num{1}                & $x_3, x_6$ & $2.5$       & $5$         \\
			Ellipse radius \num{2}                & $x_4, x_7$ & $2.5$       & $5$         \\
			Ellipse orientation (Inclusion set 1) & $x_5$      & $0$         & $\pi/2$     \\
			Ellipse orientation (Inclusion set 2) & $x_8$      & $-\pi/2$    & $\pi/2$     \\
			RVE width                             & $-$        & $15$        & $15$        \\
			RVE height                            & $x_9$      & $6$         & $15$        \\
			Initial crack location                & $w$        & $0$         & $1$         \\
			\bottomrule
		\end{tabular}
	\end{center}
\end{table}

Let $\mathbf{x} \in \mathcal{B} \subseteq \mathbb{R}^d$ be a vector of $d \in \mathbb{N}$ design variables, where $\mathcal{B}$ represents simple box constraints, with design bounds $\mathbf{x}^\text{min}, \mathbf{x}^\text{max} \in \mathbb{R}^d$.
The design variables represent the shape parameters of the inclusions represented by ellipses, cf.~Figure~\ref{fig:paper3_design_parameters}.

Let $w \in \mathcal{W} \coloneqq \{0, 0.25, 0.5, 0.75\}$ be a discrete, auxiliary variable which governs the vertical position of the initial crack relative to the microstructure, influencing the overall structural behavior.
We write the constrained optimization problem as
\begin{subequations}
	\begin{align}
    \max_{\mathbf{x} \in \mathcal{B}}  \   \Big\{ f(g(\mathbf{x},\cdot)) &\coloneqq \min_{w \in \mathcal{W}} g(\bx, w) \Big\}   ,       \\
    \text{such that} \quad                     z(\mathbf{x}) &\geq z^\text{min}, \label{eq:feasibility}
	\end{align}
  \label{eq:optimization_problem}
\end{subequations}
where $g$ is the effective fracture toughness of the material (Equation~\eqref{eq:effective_fracture_toughness}), $z$ is the clearance between inclusions and $z^\text{min}$ is a predefined minimum allowable clearance.
The clearance between the inclusion pairs are calculated using the python package \texttt{Shapely}~\cite{gillies2013shapely}.
%

We consider nine design parameters ($d=9$) related to the positions, orientations and radii of the elliptical inclusions.
The symbols and bounds of the design variables are tabulated in Table~\ref{tab:design_parameters}.
The material properties of the matrix, inclusions and interfaces are detailed in Table~\ref{tab:material_parameters}, while the phase-field model parameters are tabulated in Table~\ref{tab:algorithmic_parameters}.

\subsection{Design parametric studies}

\begin{figure}[p]
	\begin{center}
		\includegraphics[height=0.9\textheight]{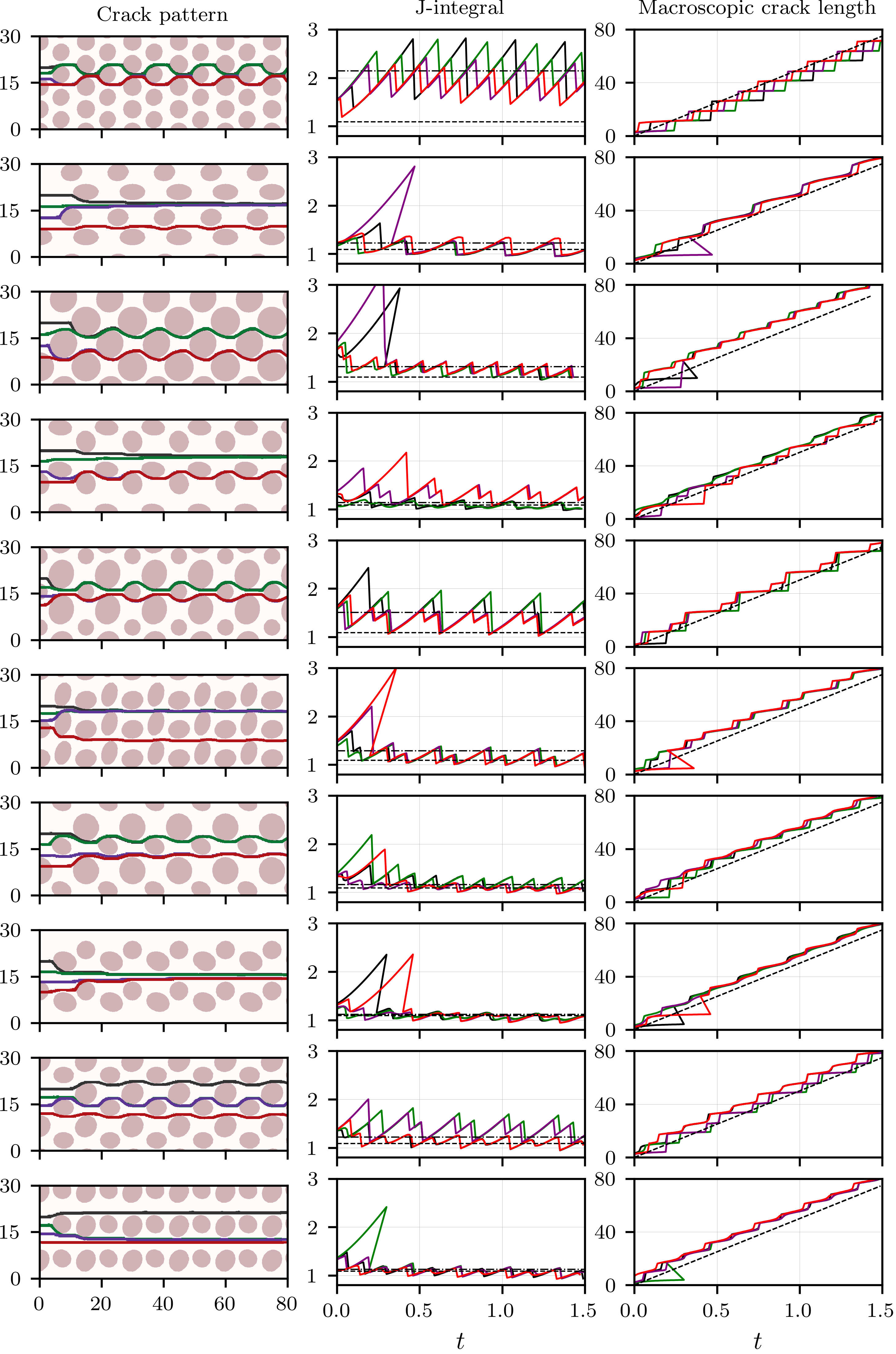}
	\end{center}
	\caption{Random designs 1-10.
	}
	\label{fig:random_designs_1}
\end{figure}

\begin{figure}[p]
	\begin{center}
		\includegraphics[height=0.9\textheight]{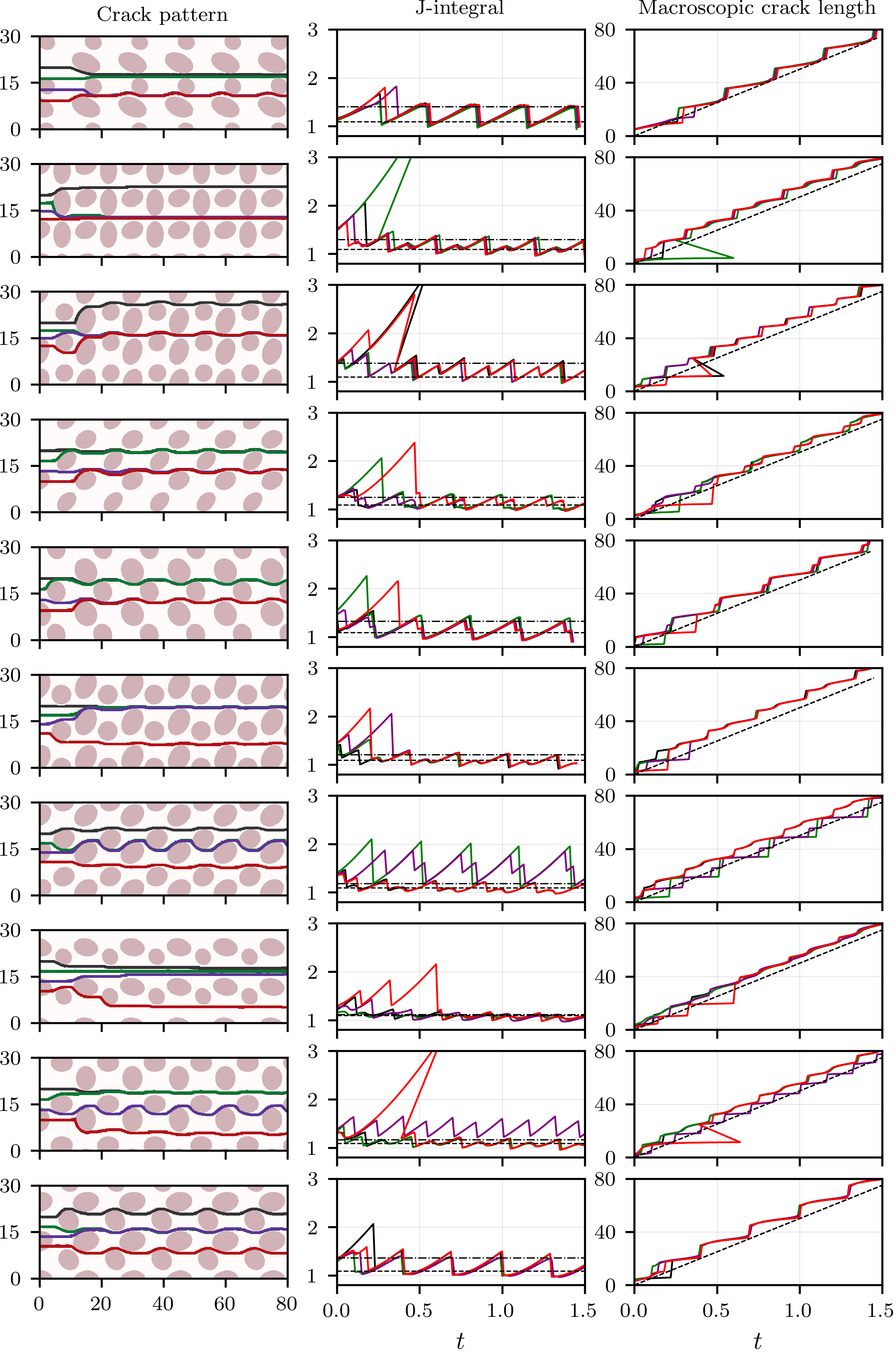}
	\end{center}
	\caption{Random designs 11-20}
	\label{fig:random_designs_2}
\end{figure}

We first conduct an initial investigation into the fracture properties of 20 randomly generated designs, all of which were filtered according to the specified feasibility criteria~\eqref{eq:feasibility}.
The results of this investigation are presented in Figure~\ref{fig:random_designs_1} and Figure~\ref{fig:random_designs_2}.
The dashed horizontal lines in the middle column of these figures correspond to the fracture toughness of the homogeneous matrix material, while the dash-dotted lines represent the effective fracture toughness of the heterogeneous material.

The fracture simulations are conducted to obtain the macroscopic crack length of \num{80}.
The effective fracture toughness of the heterogeneous material is calculated utilizing the smooth approximation of the maximum J-integral response observed for crack tip horizontal positions between \num{50} and \num{80}.

It was observed that the microstructures which allow for unobstructed crack propagation exhibit poor effective fracture toughness properties.
Conversely, when the crack is forced to propagate through obstructions caused by the inclusions, an improvement in fracture toughness is observed.
Additionally, it was determined that the position of the initial crack with respect to the inclusions has a significant influence on the fracture resistance.

Furthermore, an improvement of up to a factor of 2 was observed in the fracture toughness when compared to that of the homogeneous material.
These results suggest that the inclusion of microstructural features in the design can greatly enhance the fracture properties of the material.



\begin{figure}[ht]
	\centering
	\begin{minipage}[c]{0.5\textwidth}
	\def\svgwidth{1\columnwidth}
\begingroup%
  \makeatletter%
  \providecommand\color[2][]{%
    \errmessage{(Inkscape) Color is used for the text in Inkscape, but the package 'color.sty' is not loaded}%
    \renewcommand\color[2][]{}%
  }%
  \providecommand\transparent[1]{%
    \errmessage{(Inkscape) Transparency is used (non-zero) for the text in Inkscape, but the package 'transparent.sty' is not loaded}%
    \renewcommand\transparent[1]{}%
  }%
  \providecommand\rotatebox[2]{#2}%
  \newcommand*\fsize{\dimexpr\f@size pt\relax}%
  \newcommand*\lineheight[1]{\fontsize{\fsize}{#1\fsize}\selectfont}%
  \ifx\svgwidth\undefined%
    \setlength{\unitlength}{431.14387092bp}%
    \ifx\svgscale\undefined%
      \relax%
    \else%
      \setlength{\unitlength}{\unitlength * \real{\svgscale}}%
    \fi%
  \else%
    \setlength{\unitlength}{\svgwidth}%
  \fi%
  \global\let\svgwidth\undefined%
  \global\let\svgscale\undefined%
  \makeatother%
  \begin{picture}(1,0.7116975)%
    \lineheight{1}%
    \setlength\tabcolsep{0pt}%
    \put(0,0){\includegraphics[width=\unitlength,page=1]{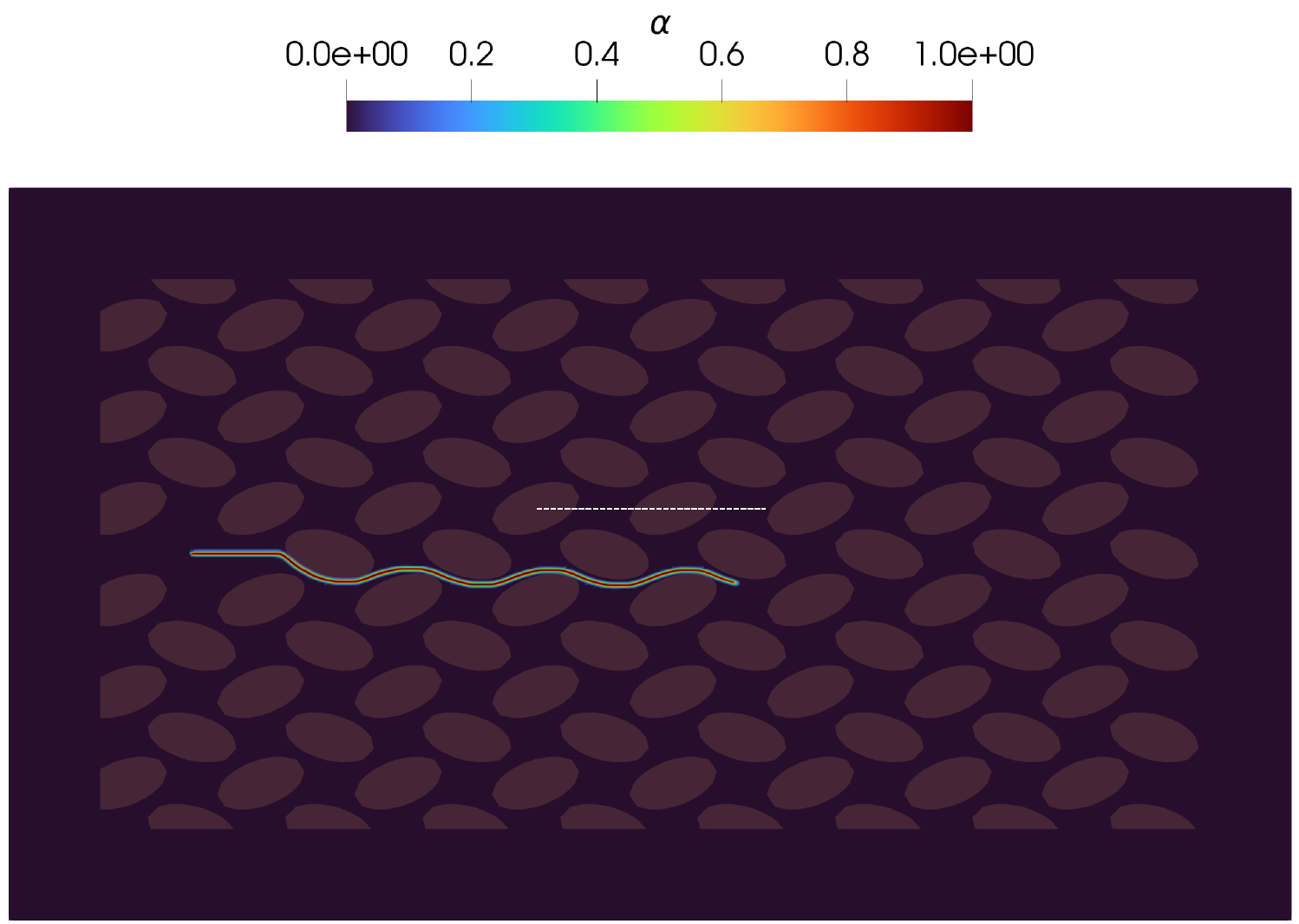}}%
    \put(0.54802434,0.3378135){\makebox(0,0)[lt]{\lineheight{1.25}\smash{\begin{tabular}[t]{l}\textcolor{white}{$\theta_1=x_5$}\end{tabular}}}}%
    \put(0,0){\includegraphics[width=\unitlength,page=2]{paper3_orientations.pdf}}%
    \put(0.81530708,0.24497446){\makebox(0,0)[lt]{\lineheight{1.25}\smash{\begin{tabular}[t]{l}\textcolor{white}{$\theta_2=x_8$}\end{tabular}}}}%
    \put(0,0){\includegraphics[width=\unitlength,page=3]{paper3_orientations.pdf}}%
  \end{picture}%
\endgroup%

	\end{minipage}
	\label{fig:paper3_orientations}
	\begin{minipage}[c]{0.39\textwidth}
		\centering
		\begin{tabular}[b]{c|c}
			\toprule
			Parameter & Value     \\
			\midrule
			\midrule
			$x_1$     & \num{0.5} \\
			$x_2$     & \num{0.5} \\
			$x_3$     & \num{5}   \\
			$x_4$     & \num{2.5} \\
			$x_6$     & \num{5}   \\
			$x_7$     & \num{2.5} \\
			$x_9$     & \num{10}  \\
			\bottomrule
		\end{tabular}
	\end{minipage}
	\caption{Geometrical setup for the two-parameter study.}
	\label{fig:geometrical_setup_two_parameter_study}
\end{figure}

\begin{figure}[ht]
	\begin{center}
		\includegraphics[width=\textwidth]{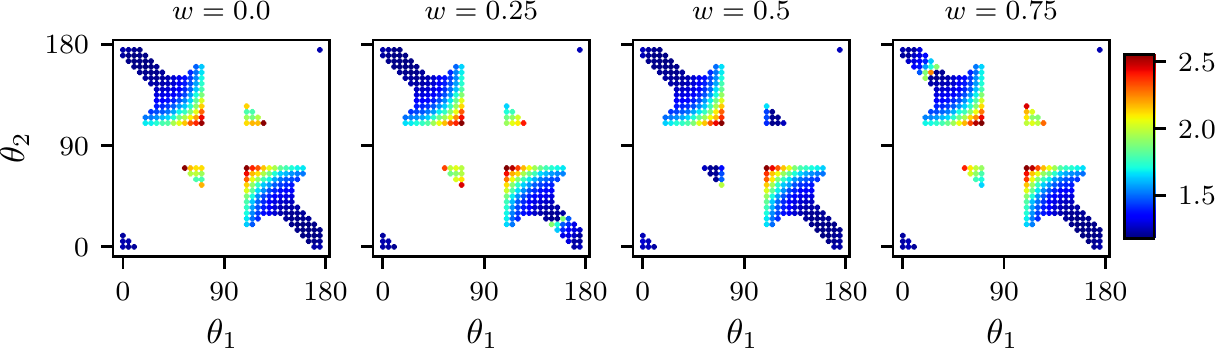}
	\end{center}
	\caption{
		Global landscapes of the objective function components with two design parameters -- $\theta_1=x_5$ and $\theta_2=x_8$ and $w=0, 0.25, 0.5, 0.75$.
		The orientations of the inclusions are measured in degrees.
		The white space in the landscape represents infeasible regions within the design space.
	}
	\label{fig:landscape_two_parameters}
\end{figure}

Next, we investigate the objective function landscape with respect to two design parameters. Specifically, we consider the orientations of the inclusions as the variable parameters, while keeping the positions and radii fixed. The geometrical setup and design parameterization for this analysis are illustrated in Figure~\ref{fig:geometrical_setup_two_parameter_study}. The landscapes of the objective function components, denoted by $g(\mathbf{x}, w_i), w_i=0, 0.25, 0.5, 0.75$, are depicted in Figure~\ref{fig:landscape_two_parameters}. It is clearly observable that the feasible design set is non-convex and disconnected, as the infeasible regions are represented by white spaces in the figures.


\subsection{Bayesian optimization}

The optimization problem~\eqref{eq:optimization_problem} is solved using the Bayesian optimization technique~\cite{ATutorialOnBFrazie2018} based on Gaussian processes~\cite{GaussianProcesRasmus2005}, which involves drawing trial points based on a belief model of the objective function $f$.
To exploit the composite nature of function $f$, employ the method proposed by \cite{BayesianOptCompositeAstudillo2019}.
The optimization procedure consists of two parts: fitting of a surrogate model to the objective function and generation of trial points based on extrema of an acquisition function subjected to inequality constraints.
In this study, we use the python package
\texttt{BoTorch}~\cite{BoTorch2019}, which is an open-source software for Bayesian optimization based on \texttt{PyTorch}~\cite{Paszke_PyTorch_An_Imperative_2019}.
For the surrogate model, \texttt{HigherOrderGP} is used and while \texttt{qUpperConfidenceBound} is employed for the acquisition function.

To obtain the next $q \in \mathbb{N}$ evaluation points $\mathbf{x}_1,\dots,\mathbf{x}_q$, the acquisition function is maximized subjected to the feasibility constraint,
which is solved using differential evolution~\cite{AConstraintHaLampin2002,Scipy10FundVirtan2020} algorithm.
For the initialization of optimization using differential evolution, $N\times q^2 \times d$ feasible design points are sampled from uniform distribution, where $N=15$ and $q=5$ were chosen for the experiments in this work.
The \texttt{differential\_evolution} implementation in \texttt{SciPy}~\cite{Scipy10FundVirtan2020} is used which allows parallel evaluation of the acquisition function as well as the constraints.

Although, Bayesian optimization algorithms are designed for searching for global optimum, the effectiveness can significantly diminish in high design dimensions.
For this reason, we also employ a successive domain reduction strategy where we define a trust region centered around the current best design.
The trust region is realized by varying the design bounds and prescribing upper and lower limits for the bounds.

The optimization was initialized using $N_0=20$ randomly generated samples.
In each iteration of the optimization process, a total of $q=5$ design points were generated. For each of these designs, $4$ simulations were performed for each value of $w$ in the set $\mathcal{W}$.
\begin{figure}[ht]
	\begin{center}
		\includegraphics[width=0.5\textwidth]{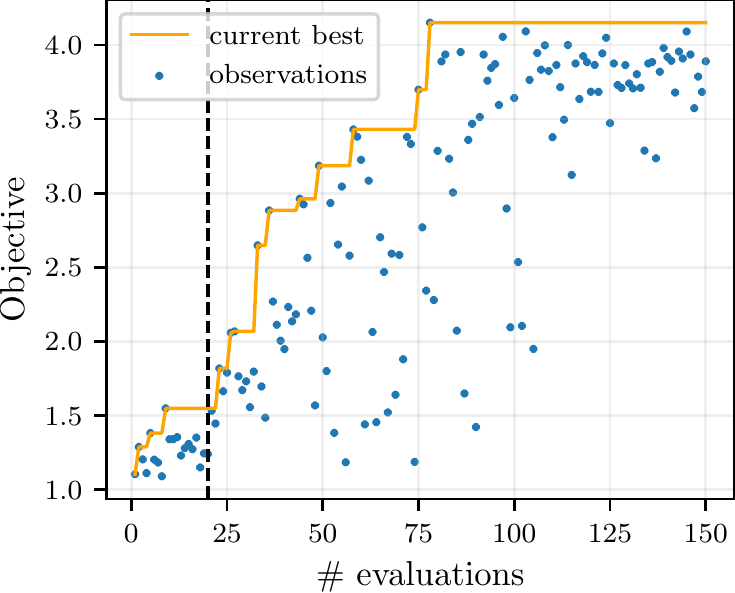}
	\end{center}
	\caption{Objective vs number of samples.
		The number of initial samples taken for fitting initial surrogate model is \num{20} and \num{5} samples are drawn in each optimization iteration.}
	\label{fig:paper3_optimization_convergence}
\end{figure}
\begin{figure}[ht]
	\begin{center}
		\includegraphics[width=\textwidth]{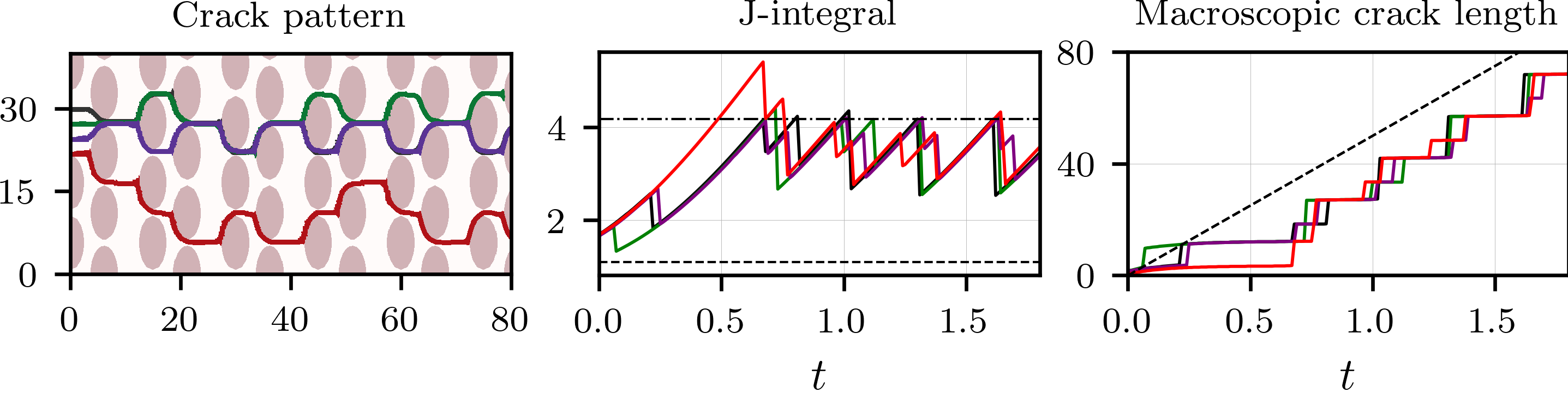}
	\end{center}
	\caption{Characteristics of the optimized design.}
	\label{fig:paper3_crack_best}
\end{figure}
Figure~\ref{fig:paper3_optimization_convergence} illustrates the convergence of the optimization process, as represented by the best objective value obtained with increasing number of evaluated design points. An improvement of $4$ times in the effective fracture toughness relative to the toughness of the homogeneous material is observed.
Figure~\ref{fig:paper3_crack_best} depicts the characteristics of the optimized design, which is found to be nearly insensitive to the location of the initial crack.

\begin{figure}[ht]
	\begin{minipage}[c]{0.49\textwidth}
		\begin{center}
			\includegraphics[width=\textwidth]{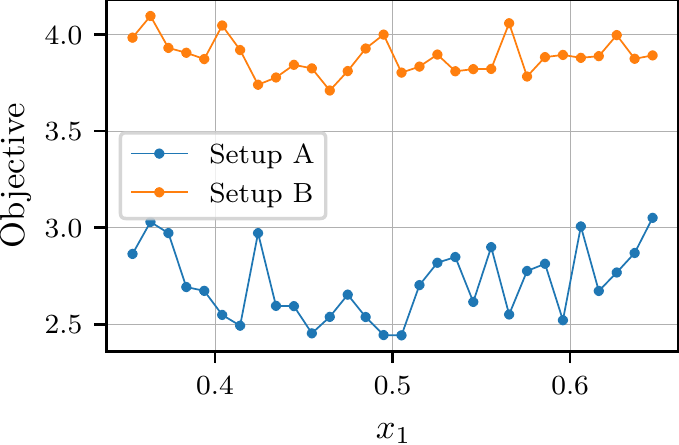}
		\end{center}
	\end{minipage}
	\begin{minipage}[c]{0.49\textwidth}
		\centering
		\begin{tabular}[b]{c|c|c}
			\toprule
			Parameter & Setup A   & Setup B   \\
			\midrule
			\midrule
			$x_2$     & \num{0.5} & \num{0.5} \\
			$x_3$     & \num{2.5} & \num{2.5} \\
			$x_4$     & \num{2.5} & \num{5}   \\
			$x_5$     & \num{0}   & \num{0}   \\
			$x_6$     & \num{2.5} & \num{2.5} \\
			$x_7$     & \num{2.5} & \num{5}   \\
			$x_8$     & \num{0}   & \num{0}   \\
			$x_9$     & \num{6}   & \num{11}  \\
			\bottomrule
		\end{tabular}
	\end{minipage}
	\caption{Objective function landscape for one design parameter for setups: A and B, as shown in Figure~\ref{fig:parametric_study_crack_visualization}.
		The noise in the objective function is due to the interaction of phase-field crack with spatial and temporal discretization.
		The truncation at the left and the right boundaries is because of the infeasibility of the design.}
	\label{fig:paper3_landscape_1d}
	\bigskip
	\centering
	\def\svgwidth{1\columnwidth}
\begingroup%
  \makeatletter%
  \providecommand\color[2][]{%
    \errmessage{(Inkscape) Color is used for the text in Inkscape, but the package 'color.sty' is not loaded}%
    \renewcommand\color[2][]{}%
  }%
  \providecommand\transparent[1]{%
    \errmessage{(Inkscape) Transparency is used (non-zero) for the text in Inkscape, but the package 'transparent.sty' is not loaded}%
    \renewcommand\transparent[1]{}%
  }%
  \providecommand\rotatebox[2]{#2}%
  \newcommand*\fsize{\dimexpr\f@size pt\relax}%
  \newcommand*\lineheight[1]{\fontsize{\fsize}{#1\fsize}\selectfont}%
  \ifx\svgwidth\undefined%
    \setlength{\unitlength}{641.71216109bp}%
    \ifx\svgscale\undefined%
      \relax%
    \else%
      \setlength{\unitlength}{\unitlength * \real{\svgscale}}%
    \fi%
  \else%
    \setlength{\unitlength}{\svgwidth}%
  \fi%
  \global\let\svgwidth\undefined%
  \global\let\svgscale\undefined%
  \makeatother%
  \begin{picture}(1,0.39446637)%
    \lineheight{1}%
    \setlength\tabcolsep{0pt}%
    \put(0,0){\includegraphics[width=\unitlength,page=1]{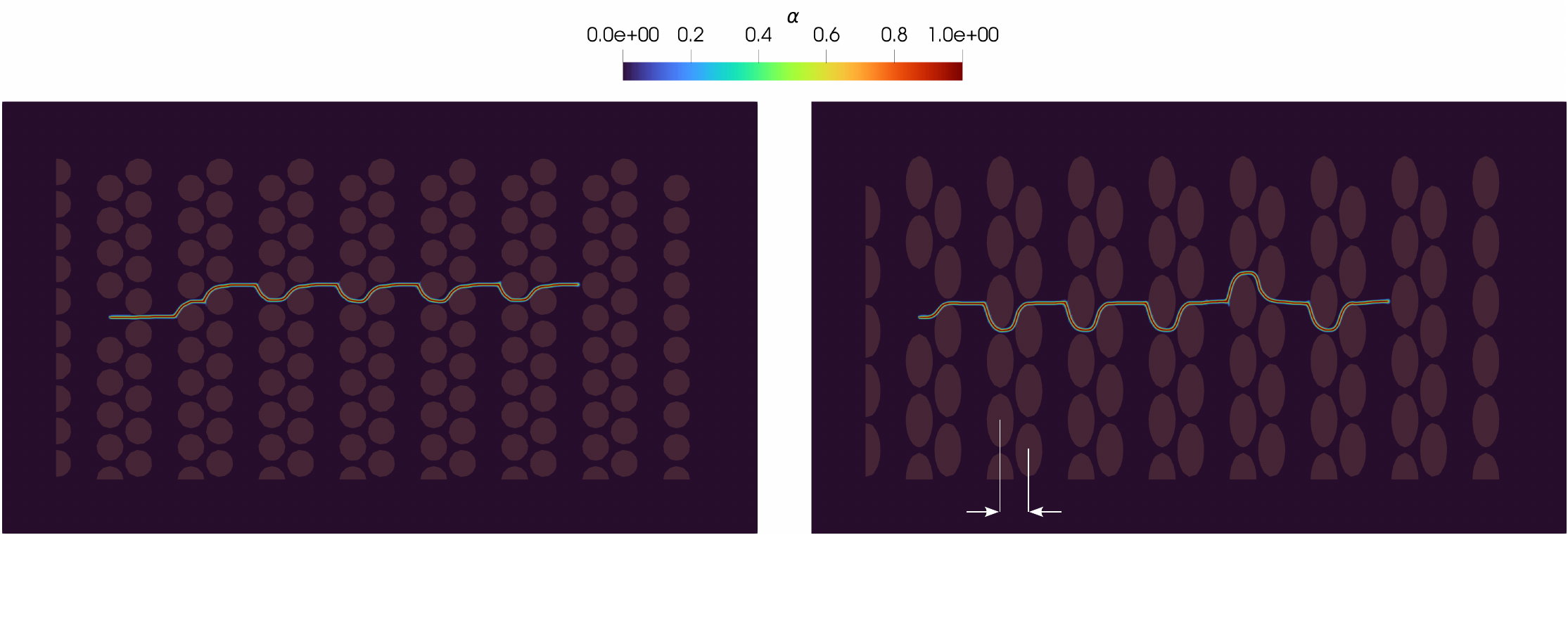}}%
    \put(0.63531347,0.03519864){\makebox(0,0)[lt]{\lineheight{1.25}\smash{\begin{tabular}[t]{l}$15x_1$\end{tabular}}}}%
    \put(0,0){\includegraphics[width=\unitlength,page=2]{parametric_study_crack_visualization.pdf}}%
    \put(0.75708533,0.03519864){\makebox(0,0)[lt]{\lineheight{1.25}\smash{\begin{tabular}[t]{l}$15$\end{tabular}}}}%
    \put(0,0){\includegraphics[width=\unitlength,page=3]{parametric_study_crack_visualization.pdf}}%
    \put(0.20970551,0.00129659){\makebox(0,0)[lt]{\lineheight{1.25}\smash{\begin{tabular}[t]{l}Setup A\end{tabular}}}}%
    \put(0.72415638,0.00129659){\makebox(0,0)[lt]{\lineheight{1.25}\smash{\begin{tabular}[t]{l}Setup B\end{tabular}}}}%
  \end{picture}%
\endgroup%

	\caption{Visualization of crack pattern with Setup A and Setup B for the parametric study with one design parameter.}
	\label{fig:parametric_study_crack_visualization}
\end{figure}

In order to establish confidence in the optimality of the design that was obtained, we conducted two additional parametric studies. In the first study, we selected circular inclusions and in the second study, we selected elongated inclusions, as shown in Figure~\ref{fig:parametric_study_crack_visualization}. The purpose of these studies was to observe the effect of horizontal spacing between the two sets of inclusion, denoted by $x_1$.

As depicted in Figure~\ref{fig:paper3_landscape_1d}, the objective function landscapes for the two cases demonstrate the presence of noise. Upon visual examination, it is apparent that there is a higher likelihood of the optimum lying towards the bounds of design variable $x_1$, as observed from the perspective of a human observer.

\section{Conclusion}\label{sec:conclusion}

In summary, this paper has presented a comprehensive framework for the optimization of microstructures to enhance fracture toughness of heterogeneous materials in two dimensions.
The proposed methodology utilizes the phase-field method to model crack propagation and employs a robust solution algorithm based on potential energy minimization and interior-point method using \texttt{IPOPT}.
The structural optimization was performed in nine design dimensions, and the sensitivity of the crack pattern to the location of the initial crack was taken into account.
The experiments conducted in this study demonstrated that the optimized designs showed improved resistance to bulk fracture.
Furthermore, design parametric studies provided additional insight into the optimality of the final design. 

\section*{Acknowledgements}
The work was funded by the
Deutsche Forschungsgemeinschaft (DFG) -- 377472739/GRK 2423/1-2019.
The authors are grateful for the administrative and advisory support from the
Competence Unit for Scientific Computing (CSC)
at
FAU Erlangen-N\"urnberg, Germany.


\printbibliography
\end{document}